\documentclass{m2an}
\usepackage{graphicx,url}

\DeclareMathOperator{\grad}{\mathbf{grad}}
\DeclareMathOperator{\dive}{\mathrm{div}}
\DeclareMathOperator{\curl}{\mathbf{curl}}

\newcommand{\ee}{\mathrm{e}}
\newcommand{\ii}{\mathrm{i}}

\newtheorem{hpthss}{Hypothesis}

\begin{document}
\title{Electromagnetic wave propagation and absorption in magnetised plasmas: variational formulations and domain decomposition}%
\thanks{This work was supported by: the Agence Nationale de la Recherche (project ``CHROME'') under contract ANR-12-BS01-0006-03; the F\'ed\'eration de Recherche Fusion par Confinement Magn\'etique--ITER; CNRS and~INRIA.}
%
\author{Aurore Back}\address{Universit\'e de Lorraine, Institut Elie Cartan de Lorraine, UMR 7502,  54506 Vand{\oe}uvre-l\`es-Nancy, France;\\
CNRS, Institut Elie Cartan de Lorraine, UMR 7502,  54506 Vand{\oe}uvre-l\`es-Nancy, France;\\ e-mail:
\url{aurore.back, takashi.hattori, simon.labrunie, jean-rodolphe.roche@univ-lorraine.fr}}
\author{Takashi Hattori}\sameaddress{1}
\author{Simon Labrunie}\sameaddress{1}
\author{Jean-Rodolphe Roche}\sameaddress{1}
\author{Pierre Bertrand}\address{Universit\'e de Lorraine, Institut Jean Lamour, UMR 7198,  54011 Nancy, France;\\
CNRS, Institut Jean Lamour, UMR 7198,  54011 Nancy, France; e-mail: \url{pierre.bertrand@univ-lorraine.fr}}
\date{January 28, 2015}
\begin{abstract}
We consider a model for the propagation and absorption of electromagnetic waves (in the time-harmonic regime) in a magnetised plasma. We present a rigorous derivation of the model and several boundary conditions modelling wave injection into the plasma. Then we propose several variational formulations, mixed and non-mixed, and prove their well-posedness thanks to a theorem by S\'ebelin \etal{} Finally, we propose a non-overlapping domain decomposition framework, show its well-posedness and equivalence with the one-domain formulation. These results appear strongly linked to the spectral properties of the plasma dielectric tensor.
\end{abstract}
\begin{resume}
Nous consid\'erons un mod\`ele de propagation et d'absorption d'ondes \'electromagn\'etiques (en r\'egime harmonique) dans un plasma magn\'etique. Nous pr\'esentons une justification rigoureuse du mod\`ele et diverses conditions aux limites mod\'elisant l'injection de l'onde dans le plasma. Puis nous proposons plusieurs formulations variationnelles, mixtes ou non, et montrons qu'elles sont bien pos\'ees gr\^ace \`a un th\'eor\`eme de S\'ebelin \etal{} Enfin, nous d\'ecrivons le principe d'une d\'ecomposition de domaine sans recouvrement, et \'etablissons le caract\`ere bien pos\'e de la formulation d\'ecompos\'ee et l'\'equivalence avec la formulation \`a un seul domaine. Ces r\'esultats paraissent intimement li\'es aux propri\'et\'es spectrales du tenseur di\'electrique du plasma. 
\end{resume}

%
\maketitle


\section{Introduction}
\label{sec-intro}

Electromagnetic wave propagation in plasmas, especially magnetised ones, is an enormous subject~\cite{Stix92}. Even in a linear framework, the equations that describe it are generally highly anisotropic and, in many practical settings, highly inhomogeneous as well. The bewildering array of phenomena and parameters involved in this modelling necessitates the derivation of simplified models tailored to the phenomenon under study, and to the theoretical or computational purpose of this study.

\medbreak

Our interest lies in the numerical simulation of the propagation of electromagnetic waves near the so-called \emph{lower hybrid} frequency in a strongly magnetised plasma. Such waves are used in tokamak technology in order to generate currents which stabilise or heat the plasma, thus bringing it closer to the conditions needed for nuclear fusion. The waves accelerate the charged particles that make up the plasma and transfer some of their energy to them through two main mechanisms: collisions between particles, which act as friction, and collisionless Landau damping. This phenomenon, caused by a resonance between electromagnetic waves and particles, is an efficient means of generating current in a magnetised plasma. Both mechanisms will be referred to as absorption. The basic physics of propagation and absorption is well understood~\cite{Stix92}. Nevertheless, efficient and robust mathematical models have to be derived in order to do reliable numerical simulations in realistic settings.

\medbreak

To perform these simulations, we have chosen to develop a finite element code which solves a suitable version of the time-harmonic Maxwell equations in a strongly magnetised plasma. It is thus a \emph{full-wave} code in the plasma community parlance, as opposed to \emph{ray-tracing} codes which solve the equations of geometrical optics. Because of their simplicity, the latter have been more popular for many years; however, it turns out that for most parameter regimes of practical interest, geometrical optics fails to hold~\cite{BrCa82}. This has renewed the interest in full-wave simulations. A full-wave code based on spectral methods in cylindrical geometry has been developed by Peysson \etal~\cite{PSL+98}. Nevertheless, generalising to real tokamak configurations, with an arbitrary cross section, requires to use the more versatile finite element method. 

\medbreak

On the other hand, full-wave computations in realistic settings are challenging because the lower hybrid wavelength is very small compared to the machine size~\cite{WBB+04}. This led us to incorporate domain decomposition capabilities in the code. This approach was preferred to, \vg, using state-of-the-art iterative methods to solve the huge linear system arising from discretisation, as it is known~\cite{ErGa12} that iterative methods perform poorly with strongly indefinite matrices such as those arising from time-harmonic equations. Thus, one has better split the computational domain into subdomains small enough to use a direct method. Another point in favour of domain decomposition is that the physical characteristics (such as density and temperature) typically vary over several orders of magnitude across a tokamak plasma. This might result in an extremely ill-conditioned linear system. However, this variation is normally continuous: another usual motivation of domain decomposition, \emph{viz.}, discontinuity in the equation coefficients, does not play any role here.

\medbreak

In this article, we will present the theoretical and mathematical foundations of our code: the derivation of the physical model, a discussion of the possible variational formulations and their well-posedness, and the domain decomposition framework. The code itself will be presented in a future publication, with a series of numerical tests. A preliminary version, without domain decomposition and differing in several respects, has been reported in~\cite{PRB+09}.

\smallbreak

The outline of the article is as follows. In~\S\ref{sec-model} we present a new, rigorous derivation of the model. We have felt this necessary, as plasma physics textbooks (such as~\cite{Stix92}) generally invoke spurious assumptions, which actually are not satisfied in real tokamak plasmas. They start by assuming that the external magnetic field and plasma characteristics are (at least approximately) homogeneous, and neglecting absorption phenomena. The latter are only discussed as an afterthought, if at all. We shall see that inhomogeneity is a non-issue, and absorption can be seamlessly integrated into the model. This is fortunate, as this phenomenon plays a crucial role in the well-posedness of the variational formulations, with or without domain decomposition. 
Actually, a simplified model related to ours is known to be ill-posed in absence of absorption~\cite{DIgW14}, which has an extremely important  consequence: the global heating effect does not vanish as absorption tends to zero. If the limiting model were well-posed, heating would be negligible when absorption is very small, as is the case in real tokamak plasmas.
The section ends with a brief discussion of the boundary conditions that can model wave injection into the plasma. 

\smallbreak

In \S\S\ref{sec-varia} and~\ref{sec-welpo} we discuss the various possible variational formulations for the injection-propagation-absorption model, and prove their well-posedness. The simplest one (which we call the \emph{plain} formulation) has been discussed in~\cite{SML+97}. We recall the results of this reference, which has not been published in a journal and thus has not reached a wide readership. Then, we show how they extend to \emph{mixed} and \emph{augmented} formulations in the sense of~\cite{Ciar05}. Mixed formulations enforce the divergence condition on the electric field, and thus control the so-called ``space charge'' phenomena. Augmented formulations allow one to use the simpler nodal finite elements~\cite[\etc]{ADH+93,AsDS96,CiZo97,Ciar05,AsSS11} instead of the edge (N\'ed\'elec) elements. Furthermore, the computed solution is continuous (like the physical solution), which avoids spurious difficulties when coupling with other solvers of computational plasma physics. We conclude this part by listing a few properties of the functional spaces that appear in the variational formulations. 

\smallbreak

In \S\ref{sec-domdec} we present the non-overlapping domain decomposition framework for our equations. Following the above discussion, we focus on the mixed augmented formulation. We prove the well-posedness of the domain-decomposed formulation and its equivalence with the initial, one-domain formulation. This parallels and generalises the work done in~\cite{AsSS11} for the time-dependent Maxwell equations in a homogeneous, isotropic, and non-absorbing medium. 
Notice, however, that our problem is considerably more difficult to solve numerically than in the latter work: the matrix of the linear system arising from the discretisation of our variational formulation is neither definite nor Hermitian, unlike that of~\cite{AsSS11}.

\section{Electromagnetic waves: model problem}
\label{sec-model}
The physical system we are interested in is a plasma or totally ionised gas, pervaded by a strong, external, static magnetic field $\boldsymbol{B}_0(\boldsymbol{x})$. (We shall always denote vector quantities by boldface letters.) Such a medium can be described as a collection of charged particles (electrons and various species of ions) which move in vacuum and create electromagnetic fields which, in turn, affect their motion. Electromagnetic fields are, thus, governed by the usual Maxwell's equations in vacuum:
\begin{eqnarray} 
\curl \boldsymbol{\mathcal{E}} = - \frac{\partial \boldsymbol{\mathcal{B}}}{\partial t} \,, &&
\curl \boldsymbol{\mathcal{B}} = \mu_0\,\boldsymbol{\mathcal{J}} +  \frac1{c^2}\, \frac{\partial \boldsymbol{\mathcal{E}}}{\partial t} \,;
\label{eq:maxrot}\\
\dive \boldsymbol{\mathcal{E}} = \varrho/\varepsilon_0, &&
\dive \boldsymbol{\mathcal{B}} = 0 .
\label{eq:maxdiv}
\end{eqnarray}
Here $\boldsymbol{\mathcal{E}}$ and $\boldsymbol{\mathcal{B}}$ denote the electric and magnetic fields; $\varrho$ and $\boldsymbol{\mathcal{J}}$  the electric charge and current densities; $\varepsilon_0$ and~$\mu_0$ the electric permittivity and magnetic permeability of vacuum, with $\varepsilon_0\,\mu_0\,c^2 =1$.

\subsection{Wave propagation equation}
The electromagnetic field is the sum of a static part and a small perturbation caused by the penetration of an electromagnetic wave. The latter is assumed to be time-harmonic. To simplify the discussion, we assume the plasma to be in mechanical and electrostatic equilibrium in the absence of the wave.
Thus, the electric and magnetic fields can be written as
\begin{equation}
\boldsymbol{\mathcal{E}}(t,\boldsymbol{x}) = \epsilon\, \Re[\boldsymbol{E}(\boldsymbol{x})\ee^{-\ii \omega t}] \quad \text{and} \quad \boldsymbol{\mathcal{B}}(t,\boldsymbol{x}) = \boldsymbol{B}_0(\boldsymbol{x}) + \epsilon\,\Re[\boldsymbol{B}(\boldsymbol{x})\ee^{-\ii \omega t}],
\label{EB-harm}
\end{equation}
where $\ii=\sqrt{-1}$, $\Re$ denotes the real part, $\epsilon \ll 1$, and $\omega > 0$ is the wave frequency. In the same way, we have
\begin{equation}
\boldsymbol{\mathcal{J}}(t,\boldsymbol{x}) = \epsilon\, \Re[\boldsymbol{J}(\boldsymbol{x})\ee^{-\ii \omega t}] \quad \text{and} \quad \varrho(t,\boldsymbol{x}) = \epsilon\, \Re[\rho(\boldsymbol{x})\ee^{-\ii \omega t}].
\label{rhoJ-harm}
\end{equation}
The static parts of $\boldsymbol{\mathcal{E}},\ \boldsymbol{\mathcal{J}}, \varrho$ are zero by the equilibrium assumption. Furthermore, the static magnetic field satisfies $\dive \boldsymbol{B}_0 = 0$ and $\curl \boldsymbol{B}_0 = \boldsymbol{0}$, as its sources are supposed to be outside the plasma. 
Plugging this ansatz in the Maxwell equations \eqref{eq:maxrot},~\eqref{eq:maxdiv}, we find:
\begin{eqnarray} 
\curl \boldsymbol{E} = \ii \omega\, \boldsymbol{B},&&
\curl \boldsymbol{B} + \ii \omega c^{-2}\, \boldsymbol{E} = \mu_{0}\, \boldsymbol{J} \,;
\label{eq:maxevolharm}\\
\dive \boldsymbol{E} = \rho/\varepsilon_0 , && \dive \boldsymbol{B} = 0.
\label{eq:maxdivharm}
\end{eqnarray}
Eliminating the variable~$\boldsymbol{B}$ between the two equations in~\eqref{eq:maxevolharm}, one finds
\begin{equation} 
\label{eqpal}
 \curl \curl \boldsymbol{E} - \tfrac{\omega^2}{c^2}  \boldsymbol{E} = \ii \omega \mu_0\, \boldsymbol{J} .
\end{equation}
We will shortly show that the medium obeys a linear, inhomogeneous and anisotropic Ohm law:
\begin{equation}
\boldsymbol{J}(\boldsymbol{x}) = \underline{\boldsymbol{\sigma}}(\boldsymbol{x})\, \boldsymbol{E}(\boldsymbol{x}) .
\label{ohm}
\end{equation}
The expression of the conductivity tensor $\underline{\boldsymbol{\sigma}}$ will be derived in the next section.
Finally, Eq.~\eqref{eqpal} becomes 
\begin{equation} 
\label{eq:principal}
\curl \curl \boldsymbol{E} - \frac{\omega^2}{c^2} \left( \underline{\boldsymbol{I}}+\frac{\ii}{\varepsilon_0 \omega}\underline{\boldsymbol{\sigma}} \right) \boldsymbol{E} = \boldsymbol{0} ,
\end{equation}
where $\underline{\boldsymbol{I}}$ is the identity matrix.

\subsection{The plasma response tensor}
As in~\cite{PSL+98,Sebe97}, the current density $\boldsymbol{J}$ in~\eqref{ohm} appears as the sum of a ``classical'' part, which can be explained by a fluid model, and a kinetic correction arising from Landau damping; both are linear in~$\boldsymbol{E}$. Let us begin with the classical part. The particles species are labelled with the subscript~$\varsigma$; the charge and mass of one particle are called $q_{\varsigma}$ and~$m_{\varsigma}$. In a first approach the plasma is assumed to be ``cold'', \ie, the thermal agitation of particles, and thus their pressure, is negligible.
Each species obeys the momentum conservation equation,%
\footnote{It can be derived by integrating in velocity the Vlasov equation, see for instance~\cite{Stix92}.}
\begin{equation}
\label{eq:transport_cold}
m_{\varsigma}\, \frac{\partial \boldsymbol{u}_{\varsigma}}{\partial t} + m_{\varsigma}\, (\boldsymbol{u}_{\varsigma} \cdot \nabla)\boldsymbol{u}_{\varsigma} - q_{\varsigma}\, (\boldsymbol{\mathcal{E}} + \boldsymbol{u}_{\varsigma} \times \boldsymbol{\mathcal{B}}) + m_{\varsigma}\, \nu_c\, \boldsymbol{u}_{\varsigma} = \boldsymbol{0}
\end{equation}
where $\boldsymbol{u}_{\varsigma}$ denotes the Eulerian fluid velocity and $\nu_c$ is the ion-electron collision frequency. (Collisions between particles of the same species do not change their bulk velocity.) The fluid velocity and the particle density $n_{\varsigma}(t,\boldsymbol{x})$ are linked to the electric charge and current densities by:
\begin{equation*}
\varrho = \sum_{\varsigma} \varrho_{\varsigma} := \sum_{\varsigma}  q_{\varsigma}\,n_{\varsigma} \,,\quad
\boldsymbol{\mathcal{J}} = \sum_{\varsigma} \boldsymbol{\mathcal{J}}_{\varsigma} := \sum_{\varsigma} q_{\varsigma}\,n_{\varsigma}\,\boldsymbol{u}_{\varsigma} \,.
\end{equation*}
Multiplying Eq.~\eqref{eq:transport_cold} by $n_{\varsigma}\,q_{\varsigma} / m_{\varsigma}$, we find
\begin{equation}
\label{eq:transport_cold_nonlinear}
\frac{\partial \boldsymbol{\mathcal{J}}_{\varsigma}}{\partial t}  + \frac1{\varrho_{\varsigma}}\, (\boldsymbol{\mathcal{J}}_{\varsigma} \cdot \nabla)\boldsymbol{\mathcal{J}}_{\varsigma} - \frac{q_{\varsigma}}{m_{\varsigma}} \left( \varrho_{\varsigma}\, \boldsymbol{\mathcal{E}} + \boldsymbol{\mathcal{J}}_{\varsigma} \times \boldsymbol{\mathcal{B}} \right) + \nu_c\, \boldsymbol{\mathcal{J}}_{\varsigma}  = \boldsymbol{0}.
\end{equation}
Then, we use the ansatz~\eqref{EB-harm}--\eqref{rhoJ-harm}. More specifically, for each species~$\varsigma$, we assume
\begin{equation*}
\boldsymbol{\mathcal{J}}_{\varsigma}(t,\boldsymbol{x}) = \epsilon\, \Re[\boldsymbol{J}_{\varsigma}(\boldsymbol{x})\ee^{-\ii \omega t}] \quad \text{and} \quad \varrho_{\varsigma}(t,\boldsymbol{x}) = q_{\varsigma}\,n_{\varsigma}^0(\boldsymbol{x}) + \epsilon\, \Re[\rho_{\varsigma}(\boldsymbol{x})\ee^{-\ii \omega t}].
\end{equation*}
The static part of~$\boldsymbol{\mathcal{J}}_{\varsigma}$ vanishes, as the plasma is at rest when $\epsilon=0$. At order~$1$ in~$\epsilon$, one can discard the term $(\boldsymbol{\mathcal{J}}_{\varsigma} \cdot \nabla)\boldsymbol{\mathcal{J}}_{\varsigma}$ altogether; in $\varrho_{\varsigma}\, \boldsymbol{\mathcal{E}}$ and~$\boldsymbol{\mathcal{J}}_{\varsigma} \times \boldsymbol{\mathcal{B}}$, only the terms in $q_{\varsigma}\,n_{\varsigma}^0\,\boldsymbol{E}$ and~$\boldsymbol{J}_{\varsigma} \times \boldsymbol{B}_0$ survive. 
Furthermore, we introduce the plasma and cyclotron frequencies for each species
\begin{equation}
\omega_{p\varsigma} := \sqrt{\frac{n_{\varsigma}^0\, q_{\varsigma}^2}{\varepsilon_0\, m_{\varsigma}}} \;,\quad  \omega_{c\varsigma} := \frac{|q_{\varsigma}|}{m_{\varsigma}} |\boldsymbol{B}_0| , 
\label{eq:omega.p.c}
\end{equation}
as well as $\delta_{\varsigma} = \mathrm{sign}(q_{\varsigma})$ and the unit vector $\boldsymbol{b} = \frac{\boldsymbol{B}_0}{|\boldsymbol{B}_0|}$.
Thus we obtain the relationship:
\begin{equation}
\ii (\omega + \ii \nu_c)\, \boldsymbol{J}_{\varsigma} + \varepsilon_0 \omega_{p\varsigma}^2\, \boldsymbol{E} + \delta_{\varsigma} \omega_{c\varsigma}\, \boldsymbol{J}_{\varsigma} \times \boldsymbol{b} = \boldsymbol{0}.
\label{eq:je1}
\end{equation}

\medbreak

At each point~$\boldsymbol{x}$, one considers an orthonormal Stix frame~\cite{Stix92} $(\boldsymbol{e}_1(\boldsymbol{x}),\boldsymbol{e}_2(\boldsymbol{x}),\boldsymbol{e}_3(\boldsymbol{x}) = \boldsymbol{b}(\boldsymbol{x}))$. For any vector field~$\boldsymbol{v}$, one denotes $\boldsymbol{v}_{\parallel} = v_3 \boldsymbol{e}_3$ and $\boldsymbol{v}_{\perp} = v_1 \boldsymbol{e}_1 + v_2 \boldsymbol{e}_2$ the components of~$\boldsymbol{v}(\boldsymbol{x})$ parallel and perpendicular to~$\boldsymbol{B}_0(\boldsymbol{x})$. 
Taking the cross product of~\eqref{eq:je1} with $\boldsymbol{b}$ on the right, we have:
\begin{equation} \label{eq:je2}
	\ii (\omega + \ii \nu_c)\, \boldsymbol{J}_{\varsigma} \times \boldsymbol{b} + \varepsilon_0 \omega_{p\varsigma}^2\, \boldsymbol{E} \times \boldsymbol{b} - \delta_{\varsigma} \omega_{c\varsigma}\, {\boldsymbol{J}_{\varsigma}}_{\perp}= \boldsymbol{0},
\end{equation}
as there holds $\boldsymbol{J}_{\perp} = \boldsymbol{b} \times (\boldsymbol{J} \times \boldsymbol{b})$.
Again, we take the cross product of~\eqref{eq:je2} with $\boldsymbol{b}$ on the left:
\begin{equation*}
\ii (\omega + \ii \nu_c)\, {\boldsymbol{J}_{\varsigma}}_{\perp} + \varepsilon_0 \omega_{p\varsigma}^2 \boldsymbol{E}_{\perp} + \delta_{\varsigma} \omega_{c\varsigma}\, \boldsymbol{J}_{\varsigma} \times \boldsymbol{b} = \boldsymbol{0},
\end{equation*}
which allows us to eliminate $\boldsymbol{J}_{\varsigma} \times \boldsymbol{b}$ in~\eqref{eq:je2}:
\begin{equation} 
\label{eq:je3}
	\left( (\omega + \ii \nu_c)^2 - \omega_{c\varsigma}^2 \right)\, {\boldsymbol{J}_{\varsigma}}_{\perp} = \ii (\omega + \ii \nu_c)\, \varepsilon_0 \omega_{p\varsigma}^2\, \boldsymbol{E}_{\perp} - \delta_{\varsigma} \omega_{c\varsigma}\, \varepsilon_0 \omega_{p\varsigma}^2\,  \boldsymbol{E} \times \boldsymbol{b}.
\end{equation}
The parallel current ${\boldsymbol{J}_{\varsigma}}_{\parallel} := (\boldsymbol{J}_{\varsigma} \cdot \boldsymbol{b})\, \boldsymbol{b}$ is obtained by taking the dot product of~\eqref{eq:je1} with~$\boldsymbol{b}$:
\begin{equation}
{\boldsymbol{J}_{\varsigma}}_{\parallel} = \frac{\ii \varepsilon_0 \omega_{p\varsigma}^2}{\omega + \ii \nu_c}\, \boldsymbol{E}_{\parallel}.
\end{equation}
Thus, the total current density $\boldsymbol{J}_{\varsigma} = {\boldsymbol{J}_{\varsigma}}_{\parallel} + {\boldsymbol{J}_{\varsigma}}_{\perp}$ of the species~$\varsigma$ is given as:
\begin{equation}
\boldsymbol{J}_{\varsigma} = \frac{\ii \varepsilon_0 \omega_{p\varsigma}^2}{\omega + \ii \nu_c} \boldsymbol{E}_{\parallel} + \frac{\ii (\omega + \ii \nu_c) \varepsilon_0 \omega_{p\varsigma}^2}{(\omega + \ii \nu_c)^2 - \omega_{c\varsigma}^2} \boldsymbol{E}_{\perp} - \frac{\delta_{\varsigma}\,\varepsilon_0 \omega_{p\varsigma}^2 \omega_{c\varsigma}}{(\omega + \ii \nu_c)^2 - \omega_{c\varsigma}^2} \boldsymbol{E} \times \boldsymbol{b}.
\end{equation}

\medbreak

Taking all species into account and setting $\alpha(\boldsymbol{x}) := \omega + \ii \nu_c(\boldsymbol{x})$, we find the expression of the ``classical'' current density:
\begin{equation}
\boldsymbol{J}_{\text{cla}} = \ii \varepsilon_0 \omega\, \underbrace{\sum_{\varsigma} \frac{ \omega_{p\varsigma}^2}{\omega \alpha }}_{=: \beta} \boldsymbol{E}_{\parallel} 
+ \ii \varepsilon_0 \omega\, \underbrace{\frac{\alpha}{\omega} \sum_{\varsigma} \frac{\omega_{p\varsigma}^{2}}{\alpha^2  - \omega_{c\varsigma}^2 }}_{=: \gamma} \boldsymbol{E}_{\perp} - \varepsilon_0 \omega\, \underbrace{\frac{1}{\omega}\, \sum_{\varsigma}\frac{\delta_{\varsigma}\,\omega_{c\varsigma}\omega^{2}_{p\varsigma}}{\alpha^{2} - \omega_{c\varsigma}^2 }}_{=: \delta} \boldsymbol{E} \times \boldsymbol{b} .
\end{equation}
In the Stix frame, we have $\boldsymbol{E} = E_1 \boldsymbol{e}_1 + E_2 \boldsymbol{e}_2 + E_3 \boldsymbol{b}$, $\boldsymbol{E}_{\parallel}  = E_3 \boldsymbol{b}$, $\boldsymbol{E} \times \boldsymbol{b} = E_2 \boldsymbol{e}_1 - E_1 \boldsymbol{e}_2$. This gives the classical part of the conductivity tensor in~\eqref{ohm}:
\begin{equation}\label{Ksigma}
\underline{\boldsymbol{\sigma}}_{\text{cla}} = \ii\varepsilon_0 \omega\,
 \begin{pmatrix}
\gamma & -\ii \delta & 0\\
\ii \delta & \gamma & 0\\
0 & 0 & \beta
\end{pmatrix} .
\end{equation}
The classical dielectric tensor is thus:
\begin{equation}\label{epsilonr}
\underline{\boldsymbol{\varepsilon}} := \boldsymbol{\underline{I}} + \frac{\ii}{\varepsilon_0 \omega}\, \underline{\boldsymbol{\sigma}}_{\text{cla}} = 
 \begin{pmatrix}
1 - \gamma & \ii \delta & 0\\
- \ii \delta & 1- \gamma & 0\\
0 & 0 & 1 - \beta
\end{pmatrix} 
:=
\begin{pmatrix}
S & -\ii D & 0 \\
\ii D & S & 0 \\
0 & 0 & P
\end{pmatrix} ,
\end{equation}
where the functions $S$, $D$ and $P$ are given by
\begin{eqnarray}
\label{coefS}
S(\boldsymbol{x}) &:=& 1 - \frac{\alpha (\boldsymbol{x})}{\omega} \sum_{\varsigma} \frac{\omega_{p\varsigma}^{2}(\boldsymbol{x})}{\alpha^2 (\boldsymbol{x}) - \omega_{c\varsigma}^2 (\boldsymbol{x})},\\ 
\label{coefD}
D(\boldsymbol{x}) &:=& \frac{1}{\omega}\, \sum_{\varsigma}\frac{\delta_{\varsigma}\,\omega_{c\varsigma}(\boldsymbol{x})\omega^{2}_{p\varsigma}(\boldsymbol{x})}{\alpha^{2}(\boldsymbol{x}) - \omega_{c\varsigma}^2 (\boldsymbol{x})},\\ 
\label{coefP}
P(\boldsymbol{x}) &:=& 1 - \sum_{\varsigma} \frac{ \omega_{p\varsigma}^2(\boldsymbol{x})}{\omega \alpha (\boldsymbol{x})}.
\end{eqnarray}

\medbreak

We proceed with the Landau damping part. As it appears~\cite{Sebe97}, only electron Landau damping in the direction parallel to~$\boldsymbol{B}_0$ plays a significant role. The ``resonant'' current generated by this effect is thus of the form:
\begin{equation}
\boldsymbol{J}_{\text{res}} (\boldsymbol{x}) = \gamma_{e}(\boldsymbol{x})\, \boldsymbol{E}_{\parallel}(\boldsymbol{x}) .
\end{equation}
The coefficient $\gamma_{e}$ is derived from a local linearisation of the Vlasov equation in the neighbourhood of the point~$\boldsymbol{x}$.
Following the classical treatment by Landau~\cite{Land46}, and assuming a Maxwellian distribution function at order~$0$ in~$\epsilon$, one finds~\cite{Sebe97}:
\begin{equation}
\gamma_e = \varepsilon_{0} \omega\, \sqrt{\frac{\pi}2}\, \frac{\omega_{pe}^2\,\omega}{k_{\parallel}^3}\, \left( \frac{m_e}{k_\mathrm{B}\,T_e} \right)^{3/2} \exp\left( -\frac{\omega^2\,m_e}{2\,k_{\parallel}^2\, k_\mathrm{B}\,T_e} \right) ,
\label{eq:gamma.e}
\end{equation}
where $T_e$ is the electron temperature (the subscript~$e$ refers to electrons), $k_\mathrm{B}$ the Boltzmann constant, and $k_{\parallel}$ is the component of the wave vector~$\boldsymbol{k}$ parallel to~$\boldsymbol{B}_0$. 

\medbreak

Adding the two contributions, $\boldsymbol{J} = \boldsymbol{J}_{\text{cla}} + \boldsymbol{J}_{\text{res}}$, we find the expression in the Stix frame of the conductivity matrix appearing in~\eqref{ohm}:
$$ \underline{\boldsymbol{\sigma}} = \underline{\boldsymbol{\sigma}}_{\text{cla}} + 
\begin{pmatrix}
0 & 0 & 0 \\
0 & 0 & 0 \\
0 & 0 & \gamma_{e}
\end{pmatrix} . $$
In other words, the equation~\eqref{eq:principal} which governs electromagnetic wave propagation and absorption in the plasma can be rewritten as:
\begin{equation} 
\curl \curl \boldsymbol{E} - \tfrac{\omega^2}{c^2} \underline{\boldsymbol{K}}\, \boldsymbol{E} = \boldsymbol{0}, 
\end{equation}
with the \emph{plasma response tensor} given by:
\begin{eqnarray}\label{Kr}
\underline{\boldsymbol{K}}(\boldsymbol{x}) = \underbrace{\left( \begin{array}{ccc}
S(\boldsymbol{x}) & -\ii D(\boldsymbol{x}) & 0 \\
  \ii D(\boldsymbol{x}) & S(\boldsymbol{x}) & 0 \\
0 & 0 & P(\boldsymbol{x})
\end{array} \right)}_{\underline{\boldsymbol{\varepsilon}}(\boldsymbol{x}) = \text{classical dielectric tensor}} + \underbrace{\frac{\ii}{\varepsilon_{0} \omega} 
\begin{pmatrix}
0 & 0 & 0 \\
0 & 0 & 0 \\
0 & 0 & \gamma_{e}(\boldsymbol{x})
\end{pmatrix}  }_{\text{Landau term}}.
\end{eqnarray}
The tensor $\underline{\boldsymbol{K}}(\boldsymbol{x})$ is not Hermitian as soon as $\nu_c$ or $\gamma_e >0$.

\subsection{The injection-propagation-absorption model}
\label{sub-model}
Let $\Omega$ be a bounded open domain in $\xR^3$, which represents the plasma volume in the tokamak. {F}rom the previous subsection, we know the propagation-absorption equation: 
\begin{equation} 
\label{eq:modela}
\curl \curl \boldsymbol{E} - \tfrac{\omega^2}{c^2} \underline{\boldsymbol{K}}\, \boldsymbol{E} = \boldsymbol{0} \quad \text{in } \Omega. 
\end{equation}
The entries of the plasma response tensor~$\underline{\boldsymbol{K}}$ are given as functions of $\boldsymbol{x} \in \Omega$.
The divergence equation
\begin{equation}
\label{eq:modelb}
\dive( \underline{\boldsymbol{K}}\, \boldsymbol{E} ) = 0 \quad \text{in } \Omega ,
\end{equation}
is a direct consequence of the previous one and may appear redundant. Nevertheless, it will play an all-important role in the derivation of the mixed and augmented variational formulations.

\begin{figure}[h]
\centerline{\includegraphics[height=4cm]{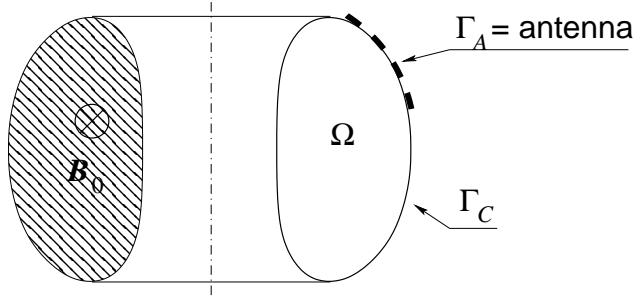}}
\caption{A cross-section of the domain $\Omega$.}
\label{fig-tok}
\end{figure}

Furthermore, various boundary conditions may be considered. Let $\Gamma$ be the boundary of~$\Omega$, and $\boldsymbol{n}$ the outgoing unitary normal vector. This boundary is made up of two parts (see Figure~\ref{fig-tok}): $\Gamma_A$~corresponds to the antenna and $\Gamma_C = \Gamma \setminus \Gamma_A$ is the remainder. Introducing the surface current $\boldsymbol{j}_A$ flowing through the antenna, the usual jump relations between media~\cite{Jack62,Boss93} give: $\boldsymbol{B}_\top = - \mu_0\, \boldsymbol{j}_A \times \boldsymbol{n}$, where $\boldsymbol{B}_\top$ denotes the  component of~$\boldsymbol{B}$ tangent to the boundary. Using the first part of~\eqref{eq:maxevolharm}, we deduce:
\begin{eqnarray}
\label{eq:eqbord01}
\curl \boldsymbol{E} \times \boldsymbol{n} = \ii \omega \mu_{0}\, \boldsymbol{j}_A \quad \textrm{on } \Gamma_{A}.
\end{eqnarray}
It appears as a Neumann (natural) condition. This modelling seems more relevant than that of~\cite{PSL+98,Sebe97,SML+97}, where $\boldsymbol{j}_A$~is treated as a ficticious volumic current in~$\Omega$.
Alternatively, one can use an essential (Dirichlet) condition:
\begin{equation}
\boldsymbol{E} \times \boldsymbol{n} = \boldsymbol{E}_A \times \boldsymbol{n} \quad \textrm{on } \Gamma_{A},
\label{eq:eqbord03}
\end{equation}
where $\boldsymbol{E}_A$ is the electric field excited at the antenna.
On the rest of the boundary, we use a perfectly conducting (homogeneous Dirichlet) boundary condition
\begin{equation}
\boldsymbol{E} \times \boldsymbol{n} = \boldsymbol{0} \quad \textrm{on } \Gamma_{C},
\label{eq:eqbord02}
\end{equation}
for the sake of simplicity.

\section{Variational formulations}
\label{sec-varia}

\subsection{Functional setting}
In the whole article, we suppose that the domain~$\Omega$ is a \emph{curved polyhedron}, \ie, a connected Lipschitz domain with piecewise smooth boundary such that, near any point of its boundary, $\Omega$~is locally $\xCinfty$-diffeomorphic to a neighbourhood of a boundary point of a polyhedron. This definition includes both smooth domains and straight polyhedra. Furthermore, the boundaries $\Gamma_A,\ \Gamma_C$ are collections of smooth faces separated by smooth edges, possibly meeting at vertices.

\medbreak

We shall use boldface letters for the functional spaces of vector fields, \vg, $\mathbf{L}^2(\Omega) = \xLtwo(\Omega)^3$. The inner product in $\mathbf{L}^2(\Omega)$ or~$\xLtwo(\Omega)$ will be denoted $(\cdot \mid \cdot)_{\Omega}$. Most unknowns and test functions are complex-valued, so this inner product is Hermitian. ``Duality'' products $\langle \varphi , v \rangle_{V}$ will be linear in the first variable~$\varphi$ and anti-linear in the second~$v$; the subscript~$V$ indicates the space to which the latter belongs. In this case, $\varphi \in V'$, the space of anti-linear forms on~$V$, which we call its dual for short. The subscripts~$\Omega,\ V$ may be dropped when the context is clear.

\medbreak

Let $\mathbf{H}(\curl;\Omega)$ be the usual space of square integrable vector fields with square integrable curl in $\Omega$. We introduce the ranges of the tangential trace mapping $\gamma_\top : \boldsymbol{v} \mapsto \boldsymbol{v} \times \boldsymbol{n}$ and the tangential component mapping $\pi_\top : \boldsymbol{v} \mapsto \boldsymbol{v}_\top := \boldsymbol{n} \times (\boldsymbol{v} \times \boldsymbol{n})$ from $\mathbf{H}(\curl;\Omega)$:
\begin{eqnarray}
\mathbf{TT}(\Gamma) &:=& \{ \boldsymbol{\varphi} \in \mathbf{H}^{-1/2}(\Gamma) : \exists \boldsymbol{v} \in \mathbf{H}(\curl;\Omega) ,\  \boldsymbol{\varphi} = \boldsymbol{v} \times \boldsymbol{n}_{|\Gamma} \},
\label{def-TN}\\
\mathbf{TC}(\Gamma) &:=& \{ \boldsymbol{\lambda} \in \mathbf{H}^{-1/2}(\Gamma) : \exists \boldsymbol{v} \in \mathbf{H}(\curl;\Omega) ,\  \boldsymbol{\lambda} = \boldsymbol{v}_{\top|\Gamma} \}.
\label{def-TT}
\end{eqnarray}
These spaces have been described in~\cite{BuCi01a}, where they are respectively denoted $\mathbf{H}^{-1/2}_{\parallel}(\dive_\Gamma,\Gamma) = \mathbf{TT}(\Gamma)$ and $\mathbf{H}^{-1/2}_{\perp}(\mathrm{curl}_\Gamma,\Gamma) = \mathbf{TC}(\Gamma)$. Furthermore~\cite{BuCi01b}, they are in duality with respect to the pivot space $\mathbf{L}^2_{t}(\Gamma) := \{ \boldsymbol{w} \in \mathbf{L}^2(\Gamma) :  \boldsymbol{w} \cdot \boldsymbol{n} = 0 \}$. This allows one to derive an integration by parts formula, valid for any $\boldsymbol{u},\ \boldsymbol{v} \in \mathbf{H}(\curl;\Omega)$:
\begin{equation}
(\boldsymbol{u} \mid \curl\boldsymbol{v})_{\Omega} - (\curl\boldsymbol{u} \mid \boldsymbol{v})_{\Omega}
= \left\langle \boldsymbol{u} \times \boldsymbol{n} , \boldsymbol{v}_\top \right\rangle_{\mathbf{TC}(\Gamma)} \,.
\label{Green2}
\end{equation}
Traces on a part of the boundary, \vg~$\Gamma_C$, can be defined straightforwardly: the range spaces, called $\mathbf{H}^{-1/2}_{\parallel,00}(\dive_{\Gamma_C},\Gamma_C)$ and $\mathbf{H}^{-1/2}_{\perp,00}(\mathrm{curl}_{\Gamma_C},\Gamma_C)$ in~\cite{BuCi01a}, will be denoted here $\mathbf{TT}(\Gamma_C),\ \mathbf{TC}(\Gamma_C)$. 
Introducing the space 
\begin{equation}
\mathbf{H}_0^{C}(\curl; \Omega) := \{ \boldsymbol{u} \in \mathbf{H}(\curl; \Omega) :  \boldsymbol{u} \times \boldsymbol{n}_{\mid \Gamma_{C}} = \boldsymbol{0} \},
\end{equation}
\ie, the subspace of fields satisfying the essential condition~\eqref{eq:eqbord02}, the range of the trace mappings on the rest of the boundary~$\Gamma_A$ will be denoted
\begin{eqnarray}
\widetilde{\mathbf{TT}}(\Gamma_A) &=& \{ \boldsymbol{\varphi} \in \mathbf{H}^{-1/2}(\Gamma_A) : \exists \boldsymbol{v} \in \mathbf{H}_0^{C}(\curl;\Omega) ,\  \boldsymbol{\varphi} = \boldsymbol{v} \times \boldsymbol{n}_{|\Gamma_A} \},
\label{def-TN-tilde}\\
&=&  \{ \boldsymbol{\varphi} \in \mathbf{TT}(\Gamma_A) : \text{the extension of $\boldsymbol{\varphi}$ by \textbf{0} to $\Gamma$ belongs to } \mathbf{TT}(\Gamma) \} ,
\nonumber\\
\widetilde{\mathbf{TC}}(\Gamma_A) &=& \{ \boldsymbol{\lambda} \in \mathbf{H}^{-1/2}(\Gamma_A) : \exists \boldsymbol{v} \in \mathbf{H}_0^{C}(\curl;\Omega) ,\  \boldsymbol{\lambda} = \boldsymbol{v}_{\top|\Gamma_A} \} 
\label{def-TT-tilde}\\
&=&  \{ \boldsymbol{\lambda} \in \mathbf{TC}(\Gamma_A) : \text{the extension of $\boldsymbol{\lambda}$ by \textbf{0} to $\Gamma$ belongs to } \mathbf{TC}(\Gamma) \} ,
\nonumber
\end{eqnarray}
instead of the ``learned'' notations $\mathbf{H}^{-1/2}_{\parallel}(\dive_{\Gamma_A}^0,\Gamma_A)$ and $\mathbf{H}^{-1/2}_{\perp}(\mathrm{curl}_{\Gamma_A}^0,\Gamma_A)$ of~\cite{BuCi01a}.
The spaces $\widetilde{\mathbf{TT}}(\Gamma_A)$ and~$\mathbf{TC}(\Gamma_A)$ are in duality with respect to the pivot space $\mathbf{L}^2_{t}(\Gamma_A)$, and similarly for $\widetilde{\mathbf{TC}}(\Gamma_A)$ and~$\mathbf{TT}(\Gamma_A)$.

\medbreak

Similarly, we introduce the Hilbert space:
\begin{equation}
\mathbf{H}(\dive \underline{\boldsymbol{K}}; \Omega) := \{ \boldsymbol{u} \in  \mathbf{L}^2(\Omega) : \dive(\underline{\boldsymbol{K}} \boldsymbol{u}) \in \xLtwo(\Omega) \} , 
\label{HdivK} 
\end{equation}
endowed with its canonical norm. If $\underline{\boldsymbol{K}} \in \xLinfty(\Omega; \mathcal{M}_3(\xC))$, it can be alternatively characterised as $\mathbf{H}(\dive \underline{\boldsymbol{K}}; \Omega) := \{ \boldsymbol{u} : \underline{\boldsymbol{K}} \boldsymbol{u} \in \mathbf{H}(\dive;\Omega) \}$, and the conormal trace $\underline{\boldsymbol{K}} \boldsymbol{u} \cdot \boldsymbol{n}$ of a field $\boldsymbol{u} \in \mathbf{H}(\dive \underline{\boldsymbol{K}}; \Omega)$ is defined as an element of~$\xHn{{-1/2}}(\Gamma)$. Then we have another useful integration by parts formula, valid for all $\boldsymbol{u} \in \mathbf{H}(\dive \underline{\boldsymbol{K}}; \Omega)$ and $\varphi\in \xHone(\Omega)$:
\begin{equation}
\left( \dive(\underline{\boldsymbol{K}} \boldsymbol{u}) \mid \varphi \right)_{\Omega} + \left( \underline{\boldsymbol{K}} \boldsymbol{u} \mid \grad \varphi \right)_{\Omega}
= \left\langle \underline{\boldsymbol{K}} \boldsymbol{u} \cdot \boldsymbol{n} , \varphi \right\rangle_{\xHn{{1/2}}(\Gamma)} \,.
\label{Green1}
\end{equation}
If $\varphi \in \xHone_0(\Omega)$, the above formula can be extended to~$\boldsymbol{u} \in \mathbf{L}^2(\Omega)$. In that case, $\dive(\underline{\boldsymbol{K}} \boldsymbol{u}) \in \xHn{{-1}}(\Omega)$ and
\begin{equation}
\left\langle \dive(\underline{\boldsymbol{K}} \boldsymbol{u}) , \varphi \right\rangle_{\xHone_0(\Omega)} + \left( \underline{\boldsymbol{K}} \boldsymbol{u} \mid \grad \varphi \right)_{\Omega} =  0.
\label{Green10}
\end{equation}
In both cases, the scalar product $\left( \underline{\boldsymbol{K}} \boldsymbol{u} \mid \grad \varphi \right)_{\Omega}$ can be written $\left( \boldsymbol{u} \mid \underline{\boldsymbol{K}}^* \grad \varphi \right)_{\Omega}$, with $\underline{\boldsymbol{K}}^*$ the conjugate transpose of~$\underline{\boldsymbol{K}}$.

\subsection{Non-mixed formulations}
Applying the Green formula~\eqref{Green2} to~\eqref{eq:modela} with the boundary conditions \eqref{eq:eqbord01} and~\eqref{eq:eqbord02}, the electric field appears as solution to the following variational formulation:\quad
{\emph{Find  $\boldsymbol{E} \in \mathbf{H}_0^{C}(\curl; \Omega)$ such that:}}
\begin{equation}
a(\boldsymbol{E},\boldsymbol{F}) = l(\boldsymbol{F}), \quad \forall \,\boldsymbol{F} \in \mathbf{H}_0^{C}(\curl; \Omega) 
\label{eq:FV}
\end{equation}
where the forms $a$ and $l$ are:
\begin{eqnarray}
a(\boldsymbol{u},\boldsymbol{v})  &:=& (\curl \boldsymbol{u}\mid \curl \boldsymbol{v})_{\Omega} - {\tfrac{\omega^2}{c^2}}( \underline{\boldsymbol{K}} \boldsymbol{u}\mid \boldsymbol{v})_{\Omega}, 
\label{eq:a01}\\
l(\boldsymbol{v}) &:=& \ii \omega \mu_{0}\langle \boldsymbol{j}_{A} , \boldsymbol{v}_{\top} \rangle_{\widetilde{\mathbf{TC}}(\Gamma_A)}.
\end{eqnarray}
The formulation~\eqref{eq:FV} will be called the \emph{plain} formulation. It can be regularised by adding to both sides a term related to the divergence. To this end, we introduce the spaces:
\begin{eqnarray}
\mathbf{X}(\underline{\boldsymbol{K}};\Omega) &:=& \mathbf{H}(\curl; \Omega) \cap \mathbf{H}( \dive \underline{\boldsymbol{K}};\Omega) ,
\label{X} \\
\mathbf{X}_N^{C}(\underline{\boldsymbol{K}};\Omega) &:=& \mathbf{H}_0^{C}(\curl; \Omega) \cap \mathbf{H}( \dive \underline{\boldsymbol{K}};\Omega) ,
\label{X0C}
\end{eqnarray}
endowed with their canonical norms and inner products. 
Using Eq.~\eqref{eq:modelb}, we obtain the \emph{augmented} variational formulation:\quad
{\emph{Find  $\boldsymbol{E} \in \mathbf{X}_N^{C}(\underline{\boldsymbol{K}};\Omega)$ such that:}}
\begin{equation}
a_{s}(\boldsymbol{E},\boldsymbol{F}) = l(\boldsymbol{F}), \quad \forall \,\boldsymbol{F} \in \mathbf{X}_N^{C}(\underline{\boldsymbol{K}};\Omega)
\label{eq:FVA}
\end{equation}
with the augmented sesquilinear form $a_{s}(\cdot,\cdot)$ ($s \in \xC$) defined on $\mathbf{X}(\underline{\boldsymbol{K}};\Omega)$ as
\begin{equation}
a_{s}(\boldsymbol{u},\boldsymbol{v}) := a(\boldsymbol{u},\boldsymbol{v}) + s\, (\dive( \underline{\boldsymbol{K}} \boldsymbol{u})\mid \dive( \underline{\boldsymbol{K}} \boldsymbol{v}))_{\Omega}.
\label{eq:as01}
\end{equation}

\subsection{Mixed formulations}
Alternatively, the divergence condition~\eqref{eq:modelb} can be considered as a constraint. Starting from the plain formulation~\eqref{eq:FV}, we introduce a Lagrangian multiplier $p \in  \xHone_0(\Omega)$ to dualise this constraint (\cf~\eqref{Green10}) and we  obtain a \emph{mixed unaugmented} formulation which writes:\quad 
\emph{Find $(\boldsymbol{E},p) \in \mathbf{H}_0^{C}(\curl; \Omega) \times \xHone_0 (\Omega)$  such that}
\begin{eqnarray} 
\label{eq:MUVFb1}
a( \boldsymbol{E}, \boldsymbol{F}) + \overline{\beta( \boldsymbol{F},p)}  &=& l( \boldsymbol{F}),\quad
\forall \boldsymbol{F} \in \mathbf{H}_0^{C}(\curl; \Omega), \\ 
\label{eq:MUVFb2}
\beta(\boldsymbol{E},q) &=& 0,\quad
\forall q \in \xHone_0(\Omega),
\end{eqnarray}
where the form $\beta$ is defined on $\mathbf{L}^2(\Omega) \times \xHone(\Omega)$ as:
\begin{equation}
\beta(\boldsymbol{v},q) := - ( \underline{\boldsymbol{K}} \boldsymbol{v} \mid \grad q)_{\Omega}.  
\label{eq;as03}
\end{equation}

\medbreak

If we start from the augmented formulation~\eqref{eq:FVA} instead, we introduce a Lagrangian multiplier $p \in  \xLtwo (\Omega)$  and arrive at the \emph{mixed augmented} variational formulation:\quad 
\emph{Find $(\boldsymbol{E},p) \in \mathbf{X}_N^{C}(\underline{\boldsymbol{K}};\Omega) \times \xLtwo (\Omega)$  such that}
\begin{eqnarray} 
\label{eq:MAVFb1}
a_{s}( \boldsymbol{E}, \boldsymbol{F}) + \overline{b( \boldsymbol{F},p)}  &=& l( \boldsymbol{F}),\quad
{\forall \boldsymbol{F} \in \mathbf{X}_N^{C}(\underline{\boldsymbol{K}};\Omega),} \\ 
\label{eq:MAVFb2}
b(\boldsymbol{E},q) &=& 0,\quad
{\forall q \in \xLtwo(\Omega),}
\end{eqnarray}
with $b(\cdot,\cdot)$ defined on $\mathbf{H}(\dive \underline{\boldsymbol{K}}; \Omega) \times \xLtwo(\Omega)$ as:
\begin{eqnarray}
b(\boldsymbol{v},q) := (\dive( \underline{\boldsymbol{K}} \boldsymbol{v})\mid q)_{\Omega}.  \label{eq;as02}
\end{eqnarray}

\subsection{Essential boundary conditions}
\label{sub-ebc}
The variational formulations for the essential conditions \eqref{eq:eqbord03} and~\eqref{eq:eqbord02} are obtained in a similar way. 
Using test fields satisfying $\boldsymbol{F} \times \boldsymbol{n} = \boldsymbol{0}$ on~$\partial\Omega$ in the Green formula~\eqref{Green2}, one derives the plain and augmented formulations:
\emph{Find $\boldsymbol{E}\in \mathbf{H}_0^C(\curl; \Omega)$, satisfying~\eqref{eq:eqbord03}, and such that}
\begin{equation}
a( \boldsymbol{E}, \boldsymbol{F}) = 0,\quad
\forall \boldsymbol{F} \in \mathbf{H}_0(\curl; \Omega) \,;
\label{eq:FVd}
\end{equation}
\emph{Find $\boldsymbol{E} \in \mathbf{X}_N^C(\underline{\boldsymbol{K}}; \Omega)$, satisfying~\eqref{eq:eqbord03}, and such that}
\begin{equation}
a_{s}( \boldsymbol{E}, \boldsymbol{F}) = 0,\quad
\forall \boldsymbol{F} \in \mathbf{X}_N(\underline{\boldsymbol{K}};\Omega),
\label{eq:FVAd}
\end{equation}
where we have set:
\begin{equation}
\mathbf{X}_N(\underline{\boldsymbol{K}}; \Omega) := \mathbf{H}_0(\curl; \Omega) \cap \mathbf{H}( \dive \underline{\boldsymbol{K}};\Omega) .
\end{equation}
To analyse these formulations, one splits the electric field as $\boldsymbol{E} = \widetilde{\boldsymbol{E}_A} + \boldsymbol{E}^\circ$, where $\widetilde{\boldsymbol{E}_A}$ is a lifting of the boundary data and $\boldsymbol{E}^\circ$ satisfies the perfectly conducting condition on the whole boundary. Thus, one has to assume at least that $\boldsymbol{E}_A \in \widetilde{\mathbf{TC}}(\Gamma_A)$.
This is obviously sufficient for the unaugmented formulations. For the augmented formulations, it is necessary to have a lifting in~$\mathbf{X}(\underline{\boldsymbol{K}};\Omega)$. As we shall see in Remark~\ref{contXhelm} the existence of such a lifting does not entail any supplementary condition on~$\boldsymbol{E}_A$. In this case, $\boldsymbol{E}^\circ$ belongs to the space~$\mathbf{X}_N(\underline{\boldsymbol{K}}; \Omega)$.
The plain and augmented formulations satisfied by~$\boldsymbol{E}^\circ$ are respectively:\\
\emph{Find $\boldsymbol{E}^\circ \in \mathbf{H}_0(\curl; \Omega)$ such that}
\begin{equation}
a( \boldsymbol{E}^\circ, \boldsymbol{F}) = \left\langle \boldsymbol{f}, \boldsymbol{F} \right\rangle_{\mathbf{H}_0(\curl; \Omega)} := -a(\widetilde{\boldsymbol{E}_A}, \boldsymbol{F}),\quad
\forall \boldsymbol{F} \in \mathbf{H}_0(\curl; \Omega) \,;
\label{eq:FVd0}
\end{equation}
\emph{Find $\boldsymbol{E}^\circ \in \mathbf{X}_N(\underline{\boldsymbol{K}}; \Omega)$ such that}
\begin{equation}
a_{s}( \boldsymbol{E}^\circ, \boldsymbol{F}) = L_{s}( \boldsymbol{F}),\quad
\forall \boldsymbol{F} \in \mathbf{X}_N(\underline{\boldsymbol{K}};\Omega).
\label{eq:FVAd0}
\end{equation}
The mixed augmented formulations satisfied by $\boldsymbol{E}$ and~$\boldsymbol{E}^\circ$ write:\\
\emph{Find $(\boldsymbol{E},p) \in \mathbf{X}_N^C(\underline{\boldsymbol{K}};\Omega) \times \xLtwo (\Omega)$, with $\boldsymbol{E}$ satisfying~\eqref{eq:eqbord03}, and such that}
\begin{eqnarray} 
\label{eq:MAVFb1d}
a_{s}( \boldsymbol{E}, \boldsymbol{F}) + \overline{b( \boldsymbol{F},p)}  &=& 0,\quad
\forall \boldsymbol{F} \in \mathbf{X}_N(\underline{\boldsymbol{K}};\Omega),\quad \\ 
\label{eq:MAVFb2d}
b(\boldsymbol{E},q) &=& 0,\quad
\forall q \in \xLtwo(\Omega).
\end{eqnarray}
\emph{Find $(\boldsymbol{E}^\circ,p) \in \mathbf{X}_N(\underline{\boldsymbol{K}};\Omega) \times \xLtwo (\Omega)$  such that}
\begin{eqnarray} 
\label{eq:MAVFb1d0}
a_{s}( \boldsymbol{E}^\circ, \boldsymbol{F}) + \overline{b( \boldsymbol{F},p)}  &=& L_{s}( \boldsymbol{F}),\quad
\forall \boldsymbol{F} \in \mathbf{X}_N(\underline{\boldsymbol{K}};\Omega),\quad \\ 
\label{eq:MAVFb2d0}
b(\boldsymbol{E}^\circ,q) &=& \ell(q),\quad
\forall q \in \xLtwo(\Omega).
\end{eqnarray}
The right-hand sides are given by:
\begin{eqnarray}
L_{s}( \boldsymbol{v}) &:=& -a(\widetilde{\boldsymbol{E}_A}, \boldsymbol{v}) - s\, (\dive( \underline{\boldsymbol{K}} \widetilde{\boldsymbol{E}_A})\mid \dive( \underline{\boldsymbol{K}} \boldsymbol{v}))_{\Omega}
:= \left\langle \boldsymbol{f}, \boldsymbol{v} \right\rangle + s\, (g \mid \dive( \underline{\boldsymbol{K}} \boldsymbol{v}))_{\Omega} \,;\qquad 
\label{rhsd1}\\
\ell(q) &:=& - (\dive( \underline{\boldsymbol{K}} \widetilde{\boldsymbol{E}_A})\mid q)_{\Omega} := (g \mid q)_{\Omega} \,.
\label{rhsd2}
\end{eqnarray}
The reader may write the mixed unaugmented formulations as an exercise.

\medbreak

As usual when dealing with non-homogeneous essential conditions, we shall use the formulations in~$\boldsymbol{E}^\circ$, Eqs.~\eqref{eq:FVd0}, \eqref{eq:FVAd0}, and~\eqref{eq:MAVFb1d0}--\eqref{eq:MAVFb2d0} to prove well-posedness. However, in practice, we discretise the formulations in~$\boldsymbol{E}$, Eqs.~\eqref{eq:FVd}, \eqref{eq:FVAd}, and~\eqref{eq:MAVFb1d}--\eqref{eq:MAVFb2d}: both conditions \eqref{eq:eqbord03} and~\eqref{eq:eqbord02} are handled by a pseudo-elimination procedure following a local change of basis~\cite{Hatt14}.

\section{Well-posedness of the problem}
\label{sec-welpo}
In this section, we summarise the results of~\cite{SML+97} --- which deals with the plain formulation when $\gamma_e=0$ --- and show how they extend to our various formulations. We shall make the following assumption throughout the article. 
\begin{hpthss}
\label{hyp:bnd}
The real functions $\nu_c$, $\omega_{c\varsigma}$ and $\omega_{p\varsigma}$, for each species~$\varsigma$ (ions and electrons) are bounded above and below by strictly positive constants on~$\Omega$. The function $\gamma_e$ is non-negative and bounded above.
\end{hpthss}
\begin{rmrk}
\label{rmq-bnd}
The collision frequency $\nu_c$ is given by the following expression~\cite{GoRu95}, where $Z$ is the ion charge number (\ie, their charge is equal to~$Z\,|q_e|$):
\begin{equation}
\nu_c = \sqrt{\frac2\pi}\, \frac{\omega_{pe}\,\ln\Lambda}{\Lambda}, \quad\text{with:}\quad \Lambda = \frac{12\pi}Z\, n_e\, \left( \frac{\varepsilon_0\,k_\mathrm{B}\,T_e}{n_e\,q_e^2} \right)^{3/2} .
\label{eq:nu.c}
\end{equation}
A plasma is characterised by $\Lambda \gg 1$. In this framework, and recalling the expressions \eqref{eq:omega.p.c} of $\omega_{c\varsigma}$,~$\omega_{p\varsigma}$ and \eqref{eq:gamma.e} of~$\gamma_e$, one checks that Hypothesis~\ref{hyp:bnd} is satisfied provided the densities $n_\varsigma$ and the electron temperature~$T_e$ are bounded above and below by strictly positive constants on~$\Omega$. This is the case in all practical settings.
\end{rmrk}

\subsection{Spectral properties of the plasma response tensor}
It is not difficult to check that the eigenvalues of the matrix $\underline{\boldsymbol{K}}$ are 
$$\lambda_1 = S+D,\quad \lambda_2 = S-D,\quad \lambda_3 = P + \dfrac{\ii}{\varepsilon_{0} \omega}\gamma_e.$$
Furthermore, $\underline{\boldsymbol{K}}$ is a normal matrix ($\underline{\boldsymbol{K}}^*\, \underline{\boldsymbol{K}} = \underline{\boldsymbol{K}}\, \underline{\boldsymbol{K}}^*$), and its singular values are the moduli of its eigenvalues. 
Then, one deduces from~\eqref{coefS}--\eqref{coefP} the expression of the imaginary parts $(\Im \lambda_i)_{i=1\ldots3}$:
\begin{eqnarray*}
\Im \lambda_1(\boldsymbol{x}) &=& \frac{\nu_c(\boldsymbol{x})}{\omega} \sum_{\varsigma}\frac{\omega^2_{p\varsigma}(\boldsymbol{x})}{(\omega^2_{c\varsigma}(\boldsymbol{x}) - \omega^2 + \nu_c^2(\boldsymbol{x}))^2 + 4\omega^2 \nu_c^2(\boldsymbol{x})}\left[(\omega - \delta_{\varsigma}\,\omega_{c\varsigma}(\boldsymbol{x}))^2 + \nu_c^2(\boldsymbol{x}) \right],\\
\Im \lambda_2(\boldsymbol{x}) &=& \frac{\nu_c(\boldsymbol{x})}{\omega} \sum_{\varsigma}\frac{\omega^2_{p\varsigma}(\boldsymbol{x})}{(\omega^2_{c\varsigma}(\boldsymbol{x}) - \omega^2 + \nu_c^2(\boldsymbol{x}))^2 + 4\omega^2 \nu_c^2(\boldsymbol{x})}\left[(\omega + \delta_{\varsigma}\,\omega_{c\varsigma}(\boldsymbol{x}))^2 + \nu_c^2(\boldsymbol{x}) \right],\\
\Im \lambda_3(\boldsymbol{x}) &=& \frac{\nu_c(\boldsymbol{x})}{\omega(\omega^2 + \nu_c^2(\boldsymbol{x}))} \sum_{\varsigma} \omega^2_{p\varsigma}(\boldsymbol{x}) + \frac{1}{\varepsilon_0 \omega} \gamma_e(\boldsymbol{x}) .  
\end{eqnarray*}
{F}rom the above calculations, one easily infers a fundamental bound.
\begin{lmm}
\label{zetaeta}
Under Hypothesis~\ref{hyp:bnd}, there exist two constants $\eta \ge \zeta > 0$, dependent on $\omega$, such that
\begin{equation}\label{zeta}
\eta (\boldsymbol{z}^* \boldsymbol{z}) \ge
|\boldsymbol{z}^* \underline{\boldsymbol{K}}(\boldsymbol{x})\boldsymbol{z}| \geq \Im[(\boldsymbol{z}^* \underline{\boldsymbol{K}}(\boldsymbol{x})\boldsymbol{z})] \geq \zeta (\boldsymbol{z}^* \boldsymbol{z}), \quad \forall \boldsymbol{z} \in \xC^3, \quad \forall \boldsymbol{x} \in \Omega.
\end{equation}
\end{lmm}
\begin{rmrk}
When the $\omega_{c\varsigma}$ and $\omega_{p\varsigma}$, $\nu_c$ and $\gamma_e$ have typical values for tokamak plasmas, and $\omega$~is of the order of the lower hybrid frequency, one has $\Re \lambda_1 \ge 0$, while $\Re \lambda_2 \le 0$ and $\Re \lambda_3 \le 0$. 
No lower bound holds for $|\Re[(\boldsymbol{z}^* \underline{\boldsymbol{K}}(\boldsymbol{x})\boldsymbol{z})]|$.
\label{rmq-zetareal}
\end{rmrk}

\subsection{Coercivity and inf-sup condition}
\label{sub-peanut}
We recall the fundamental result of S\'ebelin \etal~\cite{SML+97}.
\begin{thrm}\label{BenSeb}
Let $V$ and $H$ be Hilbert spaces such that the embedding $V \hookrightarrow H$ is continuous. Let $a(\cdot,\cdot)$ be a sesquilinear form on $V \times V$. If there exists three strictly positive constants $\alpha,\ \lambda,\ \gamma$ such that:
\begin{enumerate}
\item the real part of $a(\cdot,\cdot)$ is G{\aa}rding-elliptic on~$V$, \ie:
\begin{eqnarray}
| \Re[a(v,v)] | \geq \alpha \|v\|^2_V - \lambda \|v\|^2_H , \quad \forall v \in V
\end{eqnarray} 
\item the imaginary part of $a(\cdot,\cdot)$ is $H$-coercive, \ie:
\begin{eqnarray}
| \Im[a(v,v)] | \geq \gamma \|v\|^2_H , \quad \forall v \in V,
\end{eqnarray} 
\end{enumerate} 
then the sesquilinear form $a$ is $V$-elliptic.
\end{thrm}
Combined with Lemma~\ref{zetaeta}, this theorem shows the well-posedness of the non-mixed formulations, by Lax--Milgram's lemma.
\begin{thrm}
\label{wellposedness}
There exists a unique solution to the plain formulations \eqref{eq:FV} and~\eqref{eq:FVd0}, hence to~\eqref{eq:FVd}. The same holds for the augmented formulations \eqref{eq:FVA} and~\eqref{eq:FVAd0} --- and thus for~\eqref{eq:FVAd} --- provided $\Re s > 0$ and $\Im s \leq 0$. 
\end{thrm}
\begin{proof}
As in~\cite{SML+97}, one uses Eq.~\eqref{zeta} to check that the form $a$ given by~\eqref{eq:a01} is continuous on $V = \mathbf{H}(\curl;\Omega)$ and satisfies the assumptions of Theorem~\Rref{BenSeb} with $H = \mathbf{L}^2(\Omega)$. Thus, it is coercive (and continuous) on~$\mathbf{H}(\curl;\Omega)$, \afortiori{} on the closed subspaces $\mathbf{H}_0(\curl;\Omega)$ and~$\mathbf{H}_0^{C}(\curl;\Omega)$. 
When $\Re s > 0$ and $\Im s \leq 0$, the same applies to the form $a_{s}$ given by~\eqref{eq:as01} on $V = \mathbf{X}(\underline{\boldsymbol{K}};\Omega)$. This form is coercive and continuous on the closed subspaces of~$\mathbf{X}(\underline{\boldsymbol{K}};\Omega)$, \vg, $\mathbf{X}_N(\underline{\boldsymbol{K}};\Omega)$ and~$\mathbf{X}_N^{C}(\underline{\boldsymbol{K}};\Omega)$.
\end{proof}

\medbreak

To prove the well-posedness of mixed formulations, we have to check an inf-sup condition. This can be done by following the lines of~\cite{CiZo97,Ciar05}.
\begin{prpstn}\label{condii2}
The sesquilinear form $\beta$ defined by~\eqref{eq;as03} satisfies an inf-sup condition on $\mathbf{H}_0(\curl;\Omega) \times \xHone_0(\Omega)$, \ie there exists $C_\beta > 0$ such that
\begin{equation}
\forall q \in \xHone_0(\Omega),\quad \sup_{\boldsymbol{v} \in \mathbf{H}_0(\curl;\Omega)} \frac{|\beta(\boldsymbol{v},q)|}{\|\boldsymbol{v}\|_{\mathbf{H}(\curl)}} \geq C_{\beta}\,  \|q\|_{\xHone} \,.
\label{infsup-unaug}
\end{equation}
\end{prpstn}
\begin{proof}
Fix $q \in \xHone_0(\Omega)$ and set $\boldsymbol{v}=\grad q \in \mathbf{H}_0(\curl;\Omega)$. Lemma~\ref{zetaeta} shows that
\begin{equation*}
|\beta(\boldsymbol{v},q)| = |(\underline{\boldsymbol{K}} \boldsymbol{v} \mid \grad q )| = |(\underline{\boldsymbol{K}} \boldsymbol{v} \mid \boldsymbol{v} )| \ge \zeta \|\boldsymbol{v}\|_{\mathbf{L}^2}^2 = \zeta \|\boldsymbol{v}\|_{\mathbf{L}^2} \|\grad q\|_{\mathbf{L}^2}.
\end{equation*}
On the other hand, $\|\boldsymbol{v}\|_{\mathbf{H}(\curl)} = \big( \|\boldsymbol{v}\|_{\mathbf{L}^2}^2 + \|\curl \boldsymbol{v}\|_{\mathbf{L}^2}^2 \big)^{1/2} = \|\boldsymbol{v}\|_{\mathbf{L}^2}$, and $\|\grad q\|_{\mathbf{L}^2} \geq C\, \|q\|_{\xHone}$ by Poincar\'e's inequality. Hence the conclusion.
\end{proof}

\medbreak

\noindent To proceed to the mixed augmented case, we state and prove a useful lemma. 
\begin{lmm}
\label{divKgrad}
For any $f \in \xHn{{-1}}(\Omega)$, the elliptic problem:\quad 
\emph{Find $\phi \in \xHone_0(\Omega)$ such that}
\begin{equation}\label{divKq} 
-\Delta_{\underline{\boldsymbol{K}}}\phi := -\dive (\underline{\boldsymbol{K}} \grad \phi) = f 
\end{equation} 
admits a unique solution, which satisfies $| \psi |_{\xHone} \le C\, \| f \|_{\xHn{{-1}}} $ for some constant~$C$.
\end{lmm}
\begin{proof}
Using~\eqref{Green10} with $\boldsymbol{u} = \grad\phi$, the variational formulation of~\eqref{divKq} writes:
\begin{equation}
\mathfrak{a}(\phi,\psi) := (\underline{\boldsymbol{K}} \grad\phi \mid \grad\psi)_{\Omega} = \langle f, \psi \rangle_{\xHone_0(\Omega)},\quad \forall \psi\in \xHone_0(\Omega).
\label{fv-divKgrad}
\end{equation}
By Eq.~\eqref{zeta}, the form $\mathfrak{a}$ satisfies
$$\eta\, | \psi |_{\xHone(\Omega)}^2 \ge | \mathfrak{a}(\psi,\psi) |  \geq \zeta\, | \psi |_{\xHone(\Omega)}^2 , \quad \forall \psi\in \xHone_0(\Omega),$$
\ie, it is continuous and coercive on~$\xHone_0(\Omega)$, and the formulation is well-posed by Lax--Milgram's lemma. 
\end{proof}
\begin{rmrk}\label{divKHgrad}
Elliptic problems with the operator $\Delta_{\underline{\boldsymbol{K}}^*}$ are also well-posed.
\end{rmrk}

\begin{prpstn}\label{condii1}
The sesquilinear form $b$ defined by~\eqref{eq;as02} satisfies an inf-sup condition on $\mathbf{X}_N(\underline{\boldsymbol{K}};\Omega) \times \xLtwo(\Omega)$, \ie there exists $C_b > 0$ such that
\begin{equation}
\label{infsup-aug}
\forall q \in \xLtwo(\Omega),\quad \sup_{\boldsymbol{v} \in \mathbf{X}_N(\underline{\boldsymbol{K}};\Omega)} \frac{|b(\boldsymbol{v},q)|}{\|\boldsymbol{v}\|_{\mathbf{X}}} \geq C_{b}\,  \|q\|_{\xLtwo} \,.
\end{equation}
\end{prpstn}
\begin{proof}
Fix $q \in \xLtwo(\Omega)$. According to Lemma~\ref{divKgrad}, there exists $\phi \in \xHone_0(\Omega)$ such that $\Delta_{\underline{\boldsymbol{K}}}\phi = q$. Setting $\boldsymbol{v}=\grad \phi$, we have $\boldsymbol{v} \in \mathbf{H}_0(\curl;\Omega)$ and $\dive (\underline{\boldsymbol{K}} \boldsymbol{v}) = q$, hence $\boldsymbol{v} \in \mathbf{H}(\dive \underline{\boldsymbol{K}};\Omega)$ and finally $\boldsymbol{v} \in \mathbf{X}_N(\underline{\boldsymbol{K}};\Omega)$. It is bounded as:
\begin{eqnarray*}
\|\boldsymbol{v}\|^2_{\mathbf{X}} &=& \|\boldsymbol{v}\|_{\mathbf{L}^2}^2 + \|\curl \boldsymbol{v}\|_{\mathbf{L}^2}^2 + \|\dive \underline{\boldsymbol{K}} \boldsymbol{v}\|_{\mathbf{L}^2}^2 \nonumber \\
&=& \|\grad \phi \|_{\mathbf{L}^2}^2 + 0 + \|\dive \underline{\boldsymbol{K}} \grad \phi \|_{\xLtwo}^2 \nonumber \\
&=& | \phi |^2_{\xHone} + \|q\|_{\xLtwo}^2 \ \leq \ (1+C^2)\|q\|_{\xLtwo}^2 \,,
\end{eqnarray*}
On the other hand
\begin{equation*}
|b(\boldsymbol{v},q)| = |(\dive \underline{\boldsymbol{K}} \boldsymbol{v} \mid q )| = |(q \mid q)| = \|q\|_{\xLtwo}^2 .
\end{equation*}
Finally
$$\dfrac{|b(\boldsymbol{v},q)|}{\|\boldsymbol{v}\|_{\mathbf{X}}} \geq \dfrac{1}{\sqrt{1+C^2}} \,\|q\|_{\xLtwo} \,,$$
which we had to prove.
\end{proof}

\medbreak

\begin{thrm}
\label{wellposedness-mixed}
There exists a unique solution to the mixed unaugmented formulation~\eqref{eq:MUVFb1}--\eqref{eq:MUVFb2}, and to its counterpart for the Dirichlet boundary condition. The same holds for the mixed augmented formulations \eqref{eq:MAVFb1}--\eqref{eq:MAVFb2} and \eqref{eq:MAVFb1d0}--\eqref{eq:MAVFb2d0}, provided $\Re s > 0$ and $\Im s \leq 0$.
Thus, the problem~\eqref{eq:MAVFb1d}--\eqref{eq:MAVFb2d} is well-posed in this case.
\end{thrm}
\begin{proof}
The forms $a$ and~$a_{s}$ are coercive, in particular, on the kernels of the forms $\beta$ and~$b$. The form~$b$ is obviously continuous, and so is~$\beta(\cdot,\cdot)$ thanks to the boundedness of the entries of~$\underline{\boldsymbol{K}}$. The inf-sup conditions \eqref{infsup-unaug} and~\eqref{infsup-aug} are exactly those needed for the Dirichlet formulations. In the Neumann case, they remain valid when replacing $\mathbf{H}_0(\curl;\Omega)$ with~$\mathbf{H}_0^{C}(\curl;\Omega)$ or $\mathbf{X}_N(\underline{\boldsymbol{K}};\Omega)$ with $\mathbf{X}_N^{C}(\underline{\boldsymbol{K}};\Omega)$, as the supremum is greater on the bigger space. We conclude by the Babu\v{s}ka--Brezzi theorem. 
\end{proof}

\medbreak

To conclude this subsection, we observe that all formulations are equivalent to one another. For instance, the unique solution to the plain formulation satisfies~\eqref{eq:modela} in~$\boldsymbol{\mathcal{D}}'(\Omega)$, hence $\dive \underline{\boldsymbol{K}} \boldsymbol{E} = 0$ and $\boldsymbol{E} \in \mathbf{X}_N^{C}(\underline{\boldsymbol{K}};\Omega)$ is solution to the augmented formulation. Similarly, $(\boldsymbol{E},0)$ is solution to both mixed formulations, thus it coincides with their respective unique solutions.

\subsection{Miscellaneous properties}
Here we collect and discuss some useful properties of our functional spaces. First, one has a Helmholtz decomposition of vector fields into gradient and ``$\underline{\boldsymbol{K}}$-solenoidal'' parts.
\begin{lmm}\label{Helmvar}
For any $\boldsymbol{u} \in \mathbf{L}^2(\Omega)$ there exists a unique pair $(\phi, \boldsymbol{u}_{T}) \in \xHone_0(\Omega) \times \mathbf{L}^2(\Omega)$ satisfying the conditions
\begin{eqnarray}
\boldsymbol{u} = \grad \phi + \boldsymbol{u}_{T},&& \dive (\underline{\boldsymbol{K}} \boldsymbol{u}_{T})=0 \,; 
\label{utr}\\
\|\grad \phi\|_{\mathbf{L}^2} \leq C \|\boldsymbol{u}\|_{\mathbf{L}^2} \,, && \|\boldsymbol{u}_{T}\|_{\mathbf{L}^2} \leq C \|\boldsymbol{u}\|_{\mathbf{L}^2} \,.
\label{KHelm}
\end{eqnarray}
\end{lmm}
\begin{proof}
If a solution exists, then $\Delta_{\underline{\boldsymbol{K}}} \phi = \dive (\underline{\boldsymbol{K}} \boldsymbol{u})$; the latter function belongs to~$\xHn{{-1}}(\Omega)$ under Hypothesis~\ref{hyp:bnd}, with $\| \dive (\underline{\boldsymbol{K}} \boldsymbol{u}) \|_{\xHn{{-1}}} \le C\, \|\boldsymbol{u}\|_{\mathbf{L}^2}$. Lemma~\ref{divKgrad} shows the existence and uniqueness of~$\phi$ and $\boldsymbol{u}_T = \boldsymbol{u} - \grad \phi$, as well as the bounds~\eqref{KHelm}.
\end{proof}
\begin{rmrk}
\label{contXhelm}
Obviously, $\grad \phi \in \mathbf{H}_0(\curl;\Omega)$ and $\boldsymbol{u}_T \in \mathbf{H}(\dive\underline{\boldsymbol{K}};\Omega)$. As particular cases:
\begin{enumerate}
\item If $\boldsymbol{u} \in \mathbf{H}(\curl;\Omega)$, then $\boldsymbol{u}_T \in \mathbf{H}(\curl;\Omega)$ and thus $\boldsymbol{u}_T \in \mathbf{X}(\underline{\boldsymbol{K}};\Omega)$. Furthermore, $\boldsymbol{u}_T \times \boldsymbol{n} = \boldsymbol{u} \times \boldsymbol{n}$: the ranges of the mappings $\gamma_\top$ and~$\pi_\top$ from~$\mathbf{X}(\underline{\boldsymbol{K}};\Omega)$ are identical to the ranges from~$\mathbf{H}(\curl;\Omega)$, \ie, $\mathbf{TT}(\Gamma)$ and~$\mathbf{TC}(\Gamma)$.
\item If $\boldsymbol{u} \in \mathbf{X}(\underline{\boldsymbol{K}};\Omega)$, then $\grad\phi \in \mathbf{X}_N(\underline{\boldsymbol{K}};\Omega)$, and the decomposition~\eqref{utr} is continuous in $\mathbf{X}$~norm.
\item As a consequence of the two previous points, both $\boldsymbol{u}_T$ and $\grad\phi$ belong to~$\mathbf{X}_N(\underline{\boldsymbol{K}};\Omega)$ if $\boldsymbol{u}$ does.
\end{enumerate}
\end{rmrk}
\begin{rmrk}\label{KHgradsole}
Thanks to Remark~\ref{divKHgrad}, one also has a decomposition into $\underline{\boldsymbol{K}}^*$-gradient and solenoidal parts: for any $\boldsymbol{u} \in \mathbf{L}^2(\Omega)$, there is a unique pair $(\phi, \boldsymbol{u}_{T}) \in \xHone_0(\Omega) \times \mathbf{L}^2(\Omega)$ such that
\begin{eqnarray}
\boldsymbol{u} = \underline{\boldsymbol{K}}^* \grad \phi + \boldsymbol{u}_{T} \quad\text{and}\quad \dive  \boldsymbol{u}_{T}=0.
\end{eqnarray}
\end{rmrk}

\medbreak

The above results allow one to prove two powerful theorems on the space~$\mathbf{X}_N(\underline{\boldsymbol{K}};\Omega)$. They parallel the well-known results valid for scalar or Hermitian positive definite dielectric tensors. The proofs are similar to these classical cases and can be found in~\cite{Hatt14}, so we will not detail them here.
\begin{thrm}
\label{XTK_compacte}
If $\Omega$ is Lipschitz, the space $\mathbf{X}_N(\underline{\boldsymbol{K}};\Omega)$ is compactly embedded into~$\mathbf{L}^2(\Omega)$.
\end{thrm}
\begin{proof}
Follow the lines of Weber~\cite{Webe80}, using Lemmas \ref{zetaeta}, \ref{divKgrad} and~\ref{Helmvar}.
\end{proof}
\begin{rmrk}
If $\nu_c = \gamma_e = 0$, the proof breaks down: without absorption, one cannot establish a Fredholm alternative for the model of~\S\ref{sub-model} with Dirichlet boundary conditions, see also Remark~\ref{rmq-zetareal}. With Neumann boundary conditions, the embedding $\mathbf{X}_N^{C}(\underline{\boldsymbol{K}};\Omega) \hookrightarrow \mathbf{L}^2(\Omega)$ is \emph{not} compact when $\Gamma_A \ne \emptyset$, whatever the matrix field~$\underline{\boldsymbol{K}}$. Thus, all usual strategies for proving well-posedness fail in the absence of absorption. Actually, there is every reason to believe that the model is ill-posed in this case (see~\S\ref{sec-intro}).
\end{rmrk}
\begin{thrm}
\label{XNH1}
Assume that $\Omega$ has a $\xCn{{1,1}}$ boundary, and that the functions $\nu_c$, $\gamma_e$, $\omega_{c\varsigma}$ and $\omega_{p\varsigma}$ (for each species~$\varsigma$) belong to~$\xCone(\overline{\Omega})$. The space $\mathbf{X}_N(\underline{\boldsymbol{K}};\Omega)$ is algebraically and topologically included in~$\mathbf{H}^1(\Omega)$.
\end{thrm}
\begin{proof}
Following the lines of Birman--Solomyak~\cite{BiSo87}, one shows that any $\boldsymbol{u} \in \mathbf{X}_N(\underline{\boldsymbol{K}};\Omega)$ admits a decomposition
$$ \boldsymbol{u} = \boldsymbol{u}_{BS} + \grad\varphi,\quad \text{with:}\quad 
\boldsymbol{u}_{BS} \in \mathbf{H}^1_0(\Omega),\ \varphi \in \xHone_0(\Omega),\ \Delta_{\underline{\boldsymbol{K}}} \varphi \in \xLtwo(\Omega).$$
The usual elliptic theory~\cite{Horm84}, valid for the operator $-\Delta_{-\ii \underline{\boldsymbol{K}}}$ thanks to Lemma~\ref{zetaeta}, shows that $\varphi \in \xHn{2}(\Omega)$, given the smoothness of~$\Omega$ and~$\underline{\boldsymbol{K}}$. Hence $\boldsymbol{u} \in \mathbf{H}^1(\Omega)$.

\smallbreak

In other words, there holds $\mathbf{X}_N(\underline{\boldsymbol{K}};\Omega) \subset \mathbf{H}^1_N(\Omega) := \{ \boldsymbol{w} \in \mathbf{H}^1(\Omega) : \boldsymbol{w} \times \boldsymbol{n}_{\mid\Gamma} = \boldsymbol{0} \}$; the converse inclusion is obvious as $\underline{\boldsymbol{K}} \in \xCone(\overline{\Omega}; \mathcal{M}_3(\xC))$. Furthermore the embedding $\mathbf{H}^1_N(\Omega) \hookrightarrow \mathbf{X}_N(\underline{\boldsymbol{K}};\Omega)$ is continuous; thus the converse embedding $\mathbf{X}_N(\underline{\boldsymbol{K}};\Omega) \hookrightarrow \mathbf{H}^1_N(\Omega)$ is continuous by the open mapping theorem.
\end{proof}
\begin{rmrk}
\label{rmq-nodal-elts}
Under the hypotheses of the above theorem, it is thus possible to discretise straightforwardly the augmented and mixed augmented variational formulations of~\S\ref{sec-varia} with nodal (Lagrange or Taylor--Hood) finite elements.
Note that this does not apply when the boundary is not smooth and has re-entrant corners, due to the singularity of the solution~\cite{CoDa00}.
\end{rmrk}

\section{Non-overlapping domain decomposition framework}
\label{sec-domdec}
For the sake of simplicity, we assume from now essential boundary conditions, and we consider (\cf~\S\ref{sub-ebc}) the following model problem:
\begin{eqnarray} 
\curl \curl \boldsymbol{E} - {\tfrac{\omega^2}{c^2}} 
\underline{\boldsymbol{K}} \boldsymbol{E} &=& \boldsymbol{f} \quad \text{in } \Omega, \label{mprob1}\\ 
\dive(\underline{\boldsymbol{K}} \boldsymbol{E}) &=& g \quad \text{in } \Omega, \label{mprob2}\\
\boldsymbol{E} \times \boldsymbol{n} &=& \boldsymbol{0} \quad \text{on } \Gamma, \label{mprob3}
\end{eqnarray}
where the data $(\boldsymbol{f},g)$ satisfy the compatibility condition $\dive \boldsymbol{f} = - {\tfrac{\omega^2}{c^2}}\,g$. The mixed augmented variational formulation reads:\\
\emph{Find $(\boldsymbol{E},p) \in \mathbf{X}_N(\underline{\boldsymbol{K}};\Omega) \times \xLtwo (\Omega)$  such that}
\begin{eqnarray} 
\label{eq:fvma1}
a_{s}( \boldsymbol{E}, \boldsymbol{F}) + \overline{b( \boldsymbol{F},p)}  &=& L_{s}( \boldsymbol{F}) := ( \boldsymbol{f} \mid \boldsymbol{F} )_{\Omega} + s\, (g \mid \dive( \underline{\boldsymbol{K}} \boldsymbol{F}))_{\Omega} \,,\quad
\forall \boldsymbol{F} \in \mathbf{X}_N(\underline{\boldsymbol{K}};\Omega),\qquad \\ 
\label{eq:fvma2}
b(\boldsymbol{E},q) &=& \ell(q) := (g \mid q)_{\Omega},\quad
\forall q \in \xLtwo(\Omega).
\end{eqnarray}
As shown by the above notations, we have assumed $\boldsymbol{f} \in \mathbf{L}^2(\Omega)$ and $g \in \xLtwo(\Omega)$.
According to~\S\ref{sec-welpo}, this problem admits a unique solution $(\boldsymbol{E},p) \in \mathbf{X}_N(\underline{\boldsymbol{K}};\Omega) \times \xLtwo(\Omega)$, with $p=0$.

\subsection{Strong formulation}
We introduce a non-overlapping domain decomposition~\cite{AlVa97,QuVa99,Math08}:
\begin{equation}
\overline{\Omega} = \bigcup_{i=1}^{N_d} \overline{\Omega}_i\,;\qquad \Omega_i \subset \Omega,\quad i=1,\ldots,N_d\,;\qquad \Omega_i \cap \Omega_j = \emptyset \quad\text{if } i \neq j.
\label{dm}
\end{equation}
The exterior boundaries of subdomains are denoted $\Gamma_i=\Gamma \cap \partial\Omega_i,\ i=1,\ldots,N_d$, and the interfaces between them $\Sigma_{i,j}=\partial\Omega_i \cap \partial\Omega_j$. We shall write $i \bigtriangleup j$ whenever $\Sigma_{i,j}$ is a non-empty topological surface, \ie, it has a non-zero area.
To keep things simple, we assume that the $\Gamma_i$ and $\Sigma_{i,j}$ are smooth when they are topological surfaces, and then that the $\Gamma_i \cap \Sigma_{i,j}$ and $\Sigma_{i,j} \cap \Sigma_{i,k}$ are smooth curves. This is generally achieved in practice. However, the skeleton of interfaces $\Sigma=\bigcup_{i,j} \Sigma_{i,j}$ is not smooth, as there generally are (curved) dihedral angles between interfaces.

\medbreak 

The principle of domain decomposition for Maxwell's equations has been known for some time, both in the time-harmonic~\cite{DeJR92,AlVa97,QuVa99,Tose00,Math08} and time-dependent~\cite{AsDS96,AsSS11} versions. 
Consider the solution $\boldsymbol{E}$ to~\eqref{mprob1}--\eqref{mprob3}, and set $\boldsymbol{E}_i := \boldsymbol{E}_{|\Omega_i}$. Clearly, each $\boldsymbol{E}_i$ satisfies:    
\begin{eqnarray} 
\curl \curl \boldsymbol{E}_i - {\tfrac{\omega^2}{c^2}} 
\underline{\boldsymbol{K}} \boldsymbol{E}_i &=& \boldsymbol{f}_i \quad \textrm{ in } \Omega_i, \label{fdd}\\ 
\dive(\underline{\boldsymbol{K}}_i \boldsymbol{E}_i) &=& g_i \quad \textrm{in } \Omega_i, \label{fddconstraint}\\
\boldsymbol{E}_i \times \boldsymbol{n} &=& \boldsymbol{0} \quad \textrm{on } \Gamma_i, \label{fddext}
\end{eqnarray}
where $(\underline{\boldsymbol{K}}_i,\boldsymbol{f}_i,g_i)$ are the restrictions of $(\underline{\boldsymbol{K}},\boldsymbol{f},g)$ to $\Omega_i$.
In addition, we have the following interface conditions. As $\boldsymbol{E} \in \mathbf{H}(\curl;\Omega)$ satisfies~\eqref{mprob1} in the sense of~$\boldsymbol{\mathcal{D}}'(\Omega)$, there holds:
\begin{eqnarray*}
\boldsymbol{E}_i \times \boldsymbol{n}_i &=& - \boldsymbol{E}_j \times \boldsymbol{n}_j \quad\text{on } \Sigma_{i,j} \,, \\
\curl \boldsymbol{E}_i \times \boldsymbol{n}_i &=& - \curl  \boldsymbol{E}_j \times \boldsymbol{n}_j \quad\text{on } \Sigma_{i,j} \,,
\end{eqnarray*}
where  $\boldsymbol{n}_i$ is the outgoing unit normal vector to $\partial\Omega_i$. 
Similarly, the condition $\boldsymbol{E} \in \mathbf{H}( \dive \underline{ \boldsymbol{K}}, \Omega)$ or Eq.~\eqref{mprob2} imply
\begin{equation*}
\underline{ \boldsymbol{K}}_i \boldsymbol{E}_i \cdot \boldsymbol{n}_i = - \underline{ \boldsymbol{K}}_j  \boldsymbol{E}_j \cdot \boldsymbol{n}_j \quad\text{on } \Sigma_{i,j} \,.
\end{equation*}
Denoting as usual $[\boldsymbol{v}_i]_{\Sigma_{i,j}} = \boldsymbol{v}_i - \boldsymbol{v}_j$ (where $i$~is the larger index) the jump of~$\boldsymbol{v}$ across~$\Sigma_{i,j}$, the above interface conditions can be rewritten in the following way:
\begin{eqnarray}
&& [{ \boldsymbol{E} \times \boldsymbol{n}}]_{\Sigma_{i,j}} = 0,\qquad
[{ \underline{\boldsymbol{K}} \boldsymbol{E} \cdot \boldsymbol{n} }]_{\Sigma_{i,j}} = 0, \label{saut01}\\
&& [{ \curl \boldsymbol{E} \times \boldsymbol{n}}]_{\Sigma_{i,j}} = 0. \label{saut02}
\end{eqnarray}
Conversely, if the vector fields $\left( \boldsymbol{E}_i \right)_{i=1,\ldots,N_d}$ defined on~$\Omega_i$ satisfy Eqs.~\eqref{fdd}--\eqref{saut02} in the suitable sense, the field $\boldsymbol{E}$ defined on~$\Omega$ by glueing them is solution to~\eqref{mprob1}--\eqref{mprob3}.

\medbreak

\begin{prpstn}
\label{pro-saut3d}
Assume that the entries of~$\underline{\boldsymbol{K}}$ are continuous on~$\overline{\Omega}$. As in~\cite{AsSS11}, the interface conditions~\eqref{saut01} are equivalent to $[\boldsymbol{E}]_{\Sigma_{i,j}} = 0$. 
\end{prpstn}
\begin{proof}
The first condition implies $[\boldsymbol{E}]_{\Sigma_{i,j}} = \lambda_{i,j}\, \boldsymbol{n}_i$ for some scalar field $\lambda_{i,j}$ defined on~$\Sigma_{i,j}$. Denoting $\underline{\boldsymbol{K}}_{i,j}$ the value of~$\underline{\boldsymbol{K}}$ on this interface, the second part of~\eqref{saut01} then gives $\lambda_{i,j} \left[  \boldsymbol{n}_i \cdot (\underline{\boldsymbol{K}}_{i,j} \boldsymbol{n}_i) \right] = 0$. As $\boldsymbol{n}_i$ is a real vector, Eq.~\eqref{zeta} then implies $\lambda_{i,j} = 0$, \ie, $[\boldsymbol{E}]_{\Sigma_{i,j}} = 0$. The converse implication is obvious.
\end{proof}

\subsection{Variational formulation}
Let us now introduce a variational formulation for the multi-domain equations~\eqref{fdd}--\eqref{saut02}. The mathematical framework of domain decomposition for unaugmented Maxwell formulations is classical~\cite{AlVa97,DeJR92}. Roughly speaking, a vector field $\boldsymbol{u}_i \in \mathbf{H}(\curl;\Omega_i)$ iff it admits an extension $\boldsymbol{u} \in \mathbf{H}(\curl;\Omega)$: in this respect, $\mathbf{H}(\curl)$ spaces behave like the usual Sobolev spaces.
The case is less straightforward with augmented formulations, even when $\underline{\boldsymbol{K}} = \underline{\boldsymbol{I}}$ as in~\cite{AsSS11}. A field $\boldsymbol{u}_i \in \mathbf{H}(\curl,\dive;\Omega_i)$ does not necessarily admit an extension in~$\mathbf{H}(\curl,\dive;\Omega)$; if it does, it is actually of~$\mathbf{H}^1$ regularity, at least away from~$\Gamma_i$ when this boundary is not empty. A similar phenomenon occurs in our case. 
As said in the introduction, we shall focus on the mixed augmented formulation.

\smallbreak

We consider the following functional spaces associated to the domain decomposition~\eqref{dm}. They are endowed with their canonical ``broken'' norms. Conditions on the exterior boundary~$\Gamma_i$ are void if $\Gamma_i = \emptyset$.
\begin{eqnarray}
\mathbf{V}_0 &=& \{ \boldsymbol{v} \in \mathbf{L}^2(\Omega) : \forall i,\ \boldsymbol{v}_i := \boldsymbol{v}_{|\Omega_i} \in \mathbf{H}(\curl;\Omega_i) \text{ and } \boldsymbol{v}_i \times \boldsymbol{n} = 0  \text{ on } \Gamma_i \},  \\
\mathbf{W}_N^i &=& \{ \boldsymbol{v}_i \in \mathbf{H}(\curl;\Omega_i) \cap \mathbf{H}(\dive \underline{\boldsymbol{K}};\Omega_i) : \boldsymbol{v}_i \times \boldsymbol{n} = 0  \text{ on } \Gamma_i \}, \\
\mathbf{W}_N&=&\{ \boldsymbol{v} \in \mathbf{L}^2(\Omega) : \forall i,\ \boldsymbol{v}_{|\Omega_i} \in \mathbf{W}_N^i \}. 
\end{eqnarray}
Let $\boldsymbol{E} \in \mathbf{X}_N(\underline{\boldsymbol{K}};\Omega)$ be the solution to~\eqref{mprob1}--\eqref{mprob3}, and $\left( \boldsymbol{E}_i \right)_{i=1,\ldots,N_d}$ its decomposed version. Obviously, $\boldsymbol{E}_i  \in \mathbf{H}(\curl;\Omega_i) \cap \mathbf{H}(\dive \underline{\boldsymbol{K}};\Omega_i)$, and it satisfies~\eqref{fdd}--\eqref{fddext} as argued above. Applying the Green formula~\eqref{Green2} on each subdomain and using the first-order interface condition~\eqref{saut02}, we obtain the following formulation of Problem~\eqref{fdd}--\eqref{saut02}:
\begin{eqnarray}
\sum_i  a_{i,s}(\boldsymbol{E}_i,\boldsymbol{F}_i) + \overline{b_i(\boldsymbol{F}_i,p_i)} &=&  \sum_i L_{i,s}(\boldsymbol{F}_i), \quad \forall \, \boldsymbol{F} \in \mathbf{X}_N(\underline{\boldsymbol{K}};\Omega),    
\label{eq:ddfvma1}\\
\sum_i b_i(\boldsymbol{E}_i,q_i) &=& \sum_i \ell_i(q_i), \quad \forall \, q \in \xLtwo(\Omega ) ,
\label{eq:ddfvma2}\\ 
{} [\boldsymbol{E} \times \boldsymbol{n}]_{\Sigma_{i,j}} = 0, && [ \underline{\boldsymbol{K}} \boldsymbol{E} \cdot \boldsymbol{n} ]_{\Sigma_{i,j}} = 0, \quad \forall i \bigtriangleup j.
\label{eq:ddfvma3}
\end{eqnarray}
The domain-wise anti-linear and sesquilinear forms $a_{i,s},\ b_i,\ L_{i,s},\ \ell_i$ are defined as:
\begin{eqnarray}
a_{i,s}(\boldsymbol{u}_i,\boldsymbol{v}_i) &:=& (\curl \boldsymbol{u}_i\mid \curl \boldsymbol{v}_i)_{\Omega_i} + s\, (\dive( \underline{\boldsymbol{K}} \boldsymbol{u}_i)\mid \dive( \underline{\boldsymbol{K}} \boldsymbol{v}_i))_{\Omega_i}
- {\tfrac{\omega^2}{c^2}}( \underline{\boldsymbol{K}} \boldsymbol{u}_i\mid \boldsymbol{v}_i)_{\Omega_i} \,,
\label{eq:as01:dec}\\
b(\boldsymbol{v}_i,q_i) &:=& (\dive( \underline{\boldsymbol{K}} \boldsymbol{v}_i)\mid q_i)_{\Omega_i} \,,  \label{eq;as02:dec}\\
L_{i,s}( \boldsymbol{v}_i) &:=& ( \boldsymbol{f}_i \mid \boldsymbol{v}_i )_{\Omega_i} + s\, (g_i \mid \dive( \underline{\boldsymbol{K}} \boldsymbol{v}_i))_{\Omega_i} \,, 
\label{rhsd1:dec}\\
\ell_i(q_i) &:=& (g_i \mid q_i)_{\Omega_i} \,.
\label{rhsd2:dec}
\end{eqnarray}

\medbreak
 
In order to dualise the zeroth-order interface conditions~\eqref{saut01}, we introduce various spaces of traces and jumps. As a first step, let:
\begin{eqnarray*}
\mathbf{S}_{\Sigma}^{V} &:=& \{ \boldsymbol{\varphi} \in \mathbf{H}^{-1/2}(\Sigma) : \exists \boldsymbol{v} \in \mathbf{V}_0,\ \boldsymbol{\varphi} = [ \boldsymbol{v} \times \boldsymbol{n}]_{\Sigma} \}.
\end{eqnarray*}
The notation $[ \boldsymbol{v} \times \boldsymbol{n}]_{\Sigma}$ stands for the ordered collection of jumps $\left\{ [ \boldsymbol{v} \times \boldsymbol{n}]_{\Sigma_{i,j}} \right\}_{i \bigtriangleup j}$. Each jump belongs to~$\mathbf{TT}(\Sigma_{i,j})$, but in addition they have to satisfy some compatibility conditions~\cite{BuCi01a}. This motivates the following definition.
\begin{dfntn}
The space $\widetilde{\mathbf{TT}}(\Sigma_{i,j})$ is made of the fields $\boldsymbol{\varphi}_{i,j} \in \mathbf{TT}(\Sigma_{i,j})$ such that their extension $\boldsymbol{\varphi}$ by~$\boldsymbol{0}$ to~$\Sigma$ is the trace of a field in~$\mathbf{H}_0(\curl; \Omega)$: $\boldsymbol{\varphi} = \boldsymbol{v} \times \boldsymbol{n}_{|\Sigma}$.
\end{dfntn}
\begin{lmm}
\label{TNtilde-SSigmaV}
There holds:
$$ \bigoplus_{i \bigtriangleup j} \widetilde{\mathbf{TT}}(\Sigma_{i,j}) \subset \mathbf{S}_{\Sigma}^{V}. $$
\end{lmm}
\begin{proof}
Choose any interface $\Sigma_{i,j}$ and $\boldsymbol{\varphi}_{i,j} \in \widetilde{\mathbf{TT}}(\Sigma_{i,j})$. It can be lifted to a field $\boldsymbol{v}_i \in \mathbf{H}(\curl; \Omega_i)$ such that $\boldsymbol{v}_i \times \boldsymbol{n}_i = \boldsymbol{\varphi}_{i,j}$ on~$\Sigma_{i,j}$ and $\boldsymbol{v}_i \times \boldsymbol{n}_i = 0$ on~$\partial\Omega_i \setminus \Sigma_{i,j}$. Setting $\boldsymbol{v} = \boldsymbol{v}_i$ on~$\Omega_i$ and $\boldsymbol{0}$~elsewhere, we have $\boldsymbol{v} \in \mathbf{H}_0(\curl; \Omega)$ and $\boldsymbol{\varphi} = [\boldsymbol{v} \times \boldsymbol{n}]_{\Sigma}$. Repeating the process for all interfaces yields the conclusion.
\end{proof}

\medbreak

Then, one defines the space $\mathbb{S}_{\Sigma}^W \subset \xHn{{-1/2}}(\Sigma) \times \mathbf{S}_{\Sigma}^{V}$ as the range of the jump mapping:
\begin{eqnarray*}
\mathbf{W}_N & \longrightarrow & \xHn{{-1/2}}(\Sigma) \times \mathbf{S}_{\Sigma}^{V} \\
\boldsymbol{w} & \longmapsto & [[\boldsymbol{w} ]]_{\Sigma} :=
\left( [\underline{\boldsymbol{K}} \boldsymbol{w} \cdot \boldsymbol{n}]_{\Sigma}, [\boldsymbol{w} \times \boldsymbol{n}]_{\Sigma} \right).
\end{eqnarray*}
As in the case of~$\mathbf{S}_{\Sigma}^{V}$, the jump on~$\Sigma$ is defined by the collection of conormal and tangential jumps on the~$\Sigma_{i,j}$, which have to satisfy some compatibility conditions.
Under the assumptions of Proposition~\ref{pro-saut3d}, the data of jumps in this form is equivalent to that of three-dimensional ones.

\medbreak

One is led to introduce a new Lagrangian multiplier $\boldsymbol{\lambda} \in (\mathbb{S}_{\Sigma}^W)'$,
and we obtain the following variational formulation: \quad
{\emph{Find $(\boldsymbol{E},p,\boldsymbol{\lambda}) \in \mathbf{W}_N  \times \xLtwo(\Omega) \times (\mathbb{S}_{\Sigma}^W)'$ such that}}
\begin{eqnarray}  
\sum_i  \left\{ a_{i,s}(\boldsymbol{E}_i,\boldsymbol{F}_i) + \overline{b_i(\boldsymbol{F}_i,p_i)} \right\} + \overline{\langle \boldsymbol{\lambda} , [[\boldsymbol{F} ]]_\Sigma \rangle_{\mathbb{S}_{\Sigma}^W}}
 &=&  \sum_i L_{i,s}(\boldsymbol{F}_i), \quad \forall \, \boldsymbol{F} \in \mathbf{W}_N,  \label{eq:ddfvmac1}\\
\sum_i b_i(\boldsymbol{E}_i,q_i) &=& \sum_i \ell_i(\boldsymbol{F}_i), \quad \forall \, {q \in \xLtwo(\Omega)},  \label{eq:ddfvmac2}\\
\langle \boldsymbol{\mu} , [[ \boldsymbol{E} ]]_{\Sigma} \rangle_{\mathbb{S}_{\Sigma}^W} &=& 0, \quad \forall \boldsymbol{\mu} \in (\mathbb{S}_{\Sigma}^W)' .  
 \label{eq:ddfvmac3}
\end{eqnarray}
The duality products between $\mathbb{S}_{\Sigma}^W$ and its dual can be expressed as
\begin{eqnarray}
 \langle \boldsymbol{\lambda} , [[\boldsymbol{F} ]]_\Sigma \rangle_{\mathbb{S}_{\Sigma}^W} &=& \langle \lambda_{n} , [ \underline{\boldsymbol{K}} \boldsymbol{F} \cdot \boldsymbol{n}]_{\Sigma} \rangle + \langle \boldsymbol{\lambda}_{\top} , [\boldsymbol{F} \times \boldsymbol{n}]_{\Sigma} \rangle 
 \nonumber\\
&=& \sum_{i \bigtriangleup j} \left[ \left\langle \lambda_{n}^{i,j} \,,\, [ \underline{\boldsymbol{K}} \boldsymbol{F} \cdot \boldsymbol{n}]_{\Sigma_{i,j}} \right\rangle + \left\langle \boldsymbol{\lambda}_\top^{i,j} \,,\, [ \boldsymbol{F} \times \boldsymbol{n}]_{\Sigma_{i,j}} \right\rangle \right] 
\nonumber\\
&=& \sum_{i=1}^{N_d} \left[ \left\langle \lambda_{n} ,  \underline{\boldsymbol{K}} \boldsymbol{F}_i \cdot \boldsymbol{n}_i \right\rangle_{\xHn{{-1/2}}(\partial\Omega_i)} + \left\langle \boldsymbol{\lambda}_{\top} , \boldsymbol{F}_i \times \boldsymbol{n}_i \right\rangle_{\mathbf{TT}(\partial\Omega_i)} \right].\quad
\label{dual-Sigma}
\end{eqnarray}
On the first two lines, the dualities hold between the suitable spaces. Furthermore, on any interface $\Sigma_{i,j}$, the sum of the contributions of $\Omega_i$ and~$\Omega_j$ amounts to a jump, as the normals have opposite orientation: $\boldsymbol{n}_i = -\boldsymbol{n}_j$; hence the third line, where by convention $\lambda_n = 0$ on~$\Gamma_i$.

\begin{rmrk}
Under the assumptions of Proposition~\Rref{pro-saut3d}, the interface condition~\eqref{eq:ddfvmac3} is equivalent to:
\begin{equation}
\langle \boldsymbol{\mu} , [\boldsymbol{E}]_{\Sigma} \rangle_{\mathbf{S}_{\Sigma}^W} = 0 \quad \forall \boldsymbol{\mu} \in (\mathbf{S}_{\Sigma}^W)' ,
\end{equation}
where $\mathbf{S}_{\Sigma}^W$, the space of three-dimensional jumps of fields in~$\mathbf{W}_N$, is isomorphic to~$\mathbb{S}_{\Sigma}^W$.
If the matrix $\underline{\boldsymbol{K}} \in \xCone(\overline{\Omega}; \mathcal{M}_3(\xC))$ and $\partial\Omega$ is of~$\xCn{{1,1}}$ regularity, then $\mathbf{X}_N(\underline{\boldsymbol{K}}; \Omega) \subset \mathbf{H}^1(\Omega)$ by Theorem~\Rref{XNH1}, and the three-dimensional jump is defined in~$\mathbf{H}^{1/2}(\Sigma)$.
\end{rmrk}

\begin{rmrk}
The multi-domain variational formulation~\eqref{eq:ddfvmac1}--\eqref{eq:ddfvmac3} leads to a non-overlapping domain decomposition method. If we use nodal Taylor--Hood finite elements to discretise~\eqref{eq:ddfvmac1}--\eqref{eq:ddfvmac3}, as discussed in Remark~\ref{rmq-nodal-elts}, we obtain a saddle-point-like linear system where the unknowns are the nodal values of~$(\boldsymbol{E},p,\boldsymbol{\lambda})$. 
Mimicking Gauss factorisation, we obtain a new linear system, a generalised Schur complement system, where the unknowns are the nodal values of the Lagrange multiplier $\boldsymbol{\lambda}$ only.
To solve this new non-Hermitian reduced system,  we use a preconditioned GMRES iterative method. This algorithm induces at each iteration the resolution of a linear system corresponding to the discretisation of the variational formulation in each subdomain, as in~\cite{AsSS11}. Then we obtain a non-overlapping domain decomposition method at the discretised level. A preconditioned direct method is used to solve the linear system in each subdomain.
\end{rmrk}

\subsection{Well-posedness}
We now prove directly the well-posedness of the decomposed variational formulation.
\begin{thrm}
The decomposed formulation~\eqref{eq:ddfvmac1}--\eqref{eq:ddfvmac3} is well-posed in $\mathbf{W}_N  \times \xLtwo(\Omega) \times (\mathbb{S}_{\Sigma}^W)'$, thus it admits a unique solution.
\end{thrm}
\begin{proof}
The equations \eqref{eq:ddfvmac1}--\eqref{eq:ddfvmac3} can be written in the form of a mixed problem:\\ 
\emph{Find $(\boldsymbol{E},p,\boldsymbol{\lambda}) \in \mathbf{W}_N \times \xLtwo(\Omega) \times (\mathbb{S}_{\Sigma}^W)'$ such that}
\begin{eqnarray*}
\mathcal{A}_{s}(\boldsymbol{E},\boldsymbol{F})+\overline{\mathcal{B}(\boldsymbol{F};p,\boldsymbol{\lambda})} &=& \mathcal{L}_{s}(\boldsymbol{F}), \quad \forall \boldsymbol{F} \in \mathbf{W}_N ,\\
\mathcal{B}(\boldsymbol{E};q,\boldsymbol{\mu}) &=& \ell(q), \quad \forall (q,\boldsymbol{\mu}) \in \xLtwo(\Omega) \times (\mathbb{S}_{\Sigma}^W)',
\end{eqnarray*}  
with
\begin{eqnarray*}
\mathcal{A}_{s}(\boldsymbol{u},\boldsymbol{v}) &:=& \sum_i a_{i,s}(\boldsymbol{u}_i,\boldsymbol{v}_i) ,\\
\mathcal{B}(\boldsymbol{v};q,\boldsymbol{\mu}) &:=& \sum_i b_i(\boldsymbol{v}_i, q_i)  + \langle \boldsymbol{\mu} , [[\boldsymbol{v} ]]_{\Sigma} \rangle_{\mathbb{S}_{\Sigma}^W}.
\end{eqnarray*}
Let $\boldsymbol{v} \in \mathbf{W}_N$ and $(\boldsymbol{v}_i)_{i=1,\ldots,N_d}$ be its decomposed version. Applying Theorem~\Rref{wellposedness} in each~$\Omega_i$, one finds $a_{i,s}(\boldsymbol{v}_i,\boldsymbol{v}_i) \ge \nu_i\, \| \boldsymbol{v}_i \|_{\mathbf{W}_N^i}^2$; thus $\mathcal{A}_{s}(\boldsymbol{v},\boldsymbol{v}) \ge (\min_i\nu_i)\, \| \boldsymbol{v} \|_{\mathbf{W}_N}^2$. This holds in particular for $\boldsymbol{v} \in \ker \mathcal{B}$.

\medbreak

We denote $\mathbb{S}=\mathbb{S}_{\Sigma}^W$. To prove an inf-sup condition, we choose $(q,\boldsymbol{\mu}) \in \xLtwo(\Omega)\times \mathbb{S}'$ and seek $\boldsymbol{v} \in \mathbf{W}_N$ such that 
\begin{eqnarray}
\mathcal{B}(\boldsymbol{v};q,\boldsymbol{\mu}) \geq C_{\mathcal{B}}\, \|\boldsymbol{v}\|_{\mathbf{W}_N}\, \left( \|q\|_{\xLtwo}^2 + \|\boldsymbol{\mu}\|_{\mathbb{S}'}^2 \right)^{1/2}
\end{eqnarray}
with $C_{\mathcal{B}}$ independent of $q$ and $\boldsymbol{\mu}$.
To begin with, $\mathbb{S}$ is a space of traces, so its canonical norm is:
\begin{eqnarray*}
\|\boldsymbol{\varphi}\|_{\mathbb{S}} = \inf \left\{ \|\boldsymbol{w}\|_{\mathbf{W}_N} : [ \underline{\boldsymbol{K}} \boldsymbol{w} \cdot \boldsymbol{n}]_\Sigma = \varphi_n \text{ and }   [\boldsymbol{w} \times \boldsymbol{n} ]_\Sigma = \boldsymbol{\varphi}_{\top} \right\}.
\end{eqnarray*}
The definition of the dual norm writes:
\begin{eqnarray}
\|\boldsymbol{\mu}\|_{\mathbb{S}'} = \sup_{\boldsymbol{\varphi} \in \mathbb{S}} \dfrac{ \left| \langle \boldsymbol{\mu}, \boldsymbol{\varphi} \rangle_{\mathbb{S}} \right|}{\|\boldsymbol{\varphi}\|_{\mathbb{S}}}.
\label{dualnorm}
\end{eqnarray}
We introduce the decomposition $\mathbb{S}=\ker \boldsymbol{\mu} \oplus \xC \boldsymbol{\varphi_0}$, where $\boldsymbol{\varphi}_0$ verifies $\langle \boldsymbol{\mu}, \boldsymbol{\varphi}_0 \rangle = 1$ and $\boldsymbol{\varphi}_0 \perp \ker \boldsymbol{\mu}$. Using~\eqref{dualnorm}, we deduce $\|\boldsymbol{\varphi}_0 \|_{\mathbb{S}} = \frac{1}{\|\boldsymbol{\mu}\|_{\mathbb{S}'}}$. 

\medbreak

Consider the continuous anti-linear form $l_{\boldsymbol{\mu}}$ on $\mathbf{W}_N$ defined as
\begin{equation}
\langle l_{\boldsymbol{\mu}}, \boldsymbol{w} \rangle_{\mathbf{W}_N} =
 \langle \boldsymbol{\mu} , [[\boldsymbol{w} ]]_{\Sigma} \rangle_{\mathbb{S}} \,;
\label{elmu}
\end{equation}
obviously, it satisfies $\|l_{\boldsymbol{\mu}}\|_{\mathbf{W}_N'} \leq \|\boldsymbol{\mu}\|_{\mathbb{S}'}$. 
On the other hand, a standard argument shows the existence of a continuous lifting from $\mathbb{S}$ to $\mathbf{W}_N$ :
$$[R \boldsymbol{\varphi} \times \boldsymbol{n}]=\boldsymbol{\varphi}_{\top}, \quad  [ \underline{\boldsymbol{K}} R \boldsymbol{\varphi} \cdot \boldsymbol{n}]=\varphi_{n} \quad \textrm{and} \quad \|R \boldsymbol{\varphi}\|_{\mathbf{W}_N} \leq C_R\, \|\boldsymbol{\varphi}\|_{\mathbb{S}}.$$
Then, we introduce the decomposition $\mathbf{W}_N=\ker(l_{\boldsymbol{\mu}}) \oplus \xC \boldsymbol{w}_0$, with $\boldsymbol{w}_0= \alpha_0\, R \boldsymbol{\varphi}_0$ and $\alpha_0 \in \xC$. 
The element $\boldsymbol{w}_0 \in \mathbf{W}_N$ is normalised by the condition $\langle l_{\boldsymbol{\mu}}, \boldsymbol{w}_0 \rangle_{\mathbf{W}_N} = \|\boldsymbol{\mu}\|^2_{\mathbb{S}'}$.  This gives $\langle \boldsymbol{\mu}, \alpha_0 \boldsymbol{\varphi}_0 \rangle = \|\boldsymbol{\mu}\|^2_{\mathbb{S}'}$, hence $\alpha_0 = \|\boldsymbol{\mu}\|^2_{\mathbb{S}'}$, and finally $\|\boldsymbol{w}_0\|_{\mathbf{W}_N} \leq C_R\, \|\boldsymbol{\mu}\|_{\mathbb{S}'}$.

\medbreak

Next, consider $\phi \in \xHone_0(\Omega)$ solution to
\begin{equation*}
\Delta_{\underline{\boldsymbol{K}}}\phi = f \in \xLtwo(\Omega),\quad \text{with} \quad f_i=q_i - \dive (\underline{\boldsymbol{K}}\boldsymbol{w}_{0i}) \text{ in } \Omega_i \,.
\end{equation*}
This function is bounded as:
\begin{eqnarray*}
|\phi|_{\xHone(\Omega)} \leq  C_1\, (\|q\|_{\xLtwo(\Omega)} + \|\boldsymbol{w}_0\|_{\mathbf{W}_N})  \leq  C_1'\, (\|q\|_{\xLtwo(\Omega)} + \|\boldsymbol{\mu}\|_{\mathbb{S}'}) .
\end{eqnarray*}
The vector field $\boldsymbol{v} := \boldsymbol{w}_0 + \grad \phi \in \mathbf{W}_N$ satifies $\dive \underline{\boldsymbol{K}} \boldsymbol{v}_i = q_i$ in~$\Omega_i$, and is bounded as:
\begin{eqnarray*}
\|\boldsymbol{v}\|_{\mathbf{W}_N}^2 &=& \|\boldsymbol{w}_0 + \grad \phi\|^2_{\xLtwo(\Omega)} + \sum_i \Big[ \|\curl \boldsymbol{w}_0\|^2_{\xLtwo(\Omega_i)} + \|q\|^2_{\xLtwo(\Omega_i)} \Big] \\
& \leq & 2 \big[ \|\boldsymbol{w}_0\|^2_{\mathbf{W}_N} + |\phi|^2_{\xHone(\Omega)} \big] + \|q\|^2_{\xLtwo(\Omega)} 
\ \leq \ C_2 \big( \|q\|^2_{\xLtwo(\Omega)} + \|\boldsymbol{\mu}\|^2_{\mathbb{S}'} \big).
\end{eqnarray*}
But $\grad \phi \in \mathbf{X}_N(\underline{\boldsymbol{K}};\Omega)$, which implies:
\begin{eqnarray*}
[\underline{\boldsymbol{K}} \boldsymbol{v} \cdot \boldsymbol{n}]_\Sigma = [\underline{\boldsymbol{K}} \boldsymbol{w}_0 \cdot \boldsymbol{n}]_\Sigma \quad \text{and} \quad [\boldsymbol{v} \times \boldsymbol{n}]_\Sigma = [\boldsymbol{w}_0 \times \boldsymbol{n}]_\Sigma.
\end{eqnarray*}
We conclude that
\begin{eqnarray*}
\mathcal{B}(\boldsymbol{v};q,\boldsymbol{\mu}) &=& \sum_i (\dive \underline{\boldsymbol{K}} \boldsymbol{v}_i \mid q_i) + \langle \boldsymbol{\mu} , [[\boldsymbol{v}]]_\Sigma \rangle_{\mathbb{S}}, \\
&=& \sum_i \|q_i\|^2_{\xLtwo(\Omega_i)} + \underbrace{\langle \boldsymbol{\mu} , [[\boldsymbol{w}_0]]_\Sigma \rangle_{\mathbb{S}}}_{\langle l_{\boldsymbol{\mu}}, \boldsymbol{w}_0 \rangle_{\mathbf{W}_N} = \|\boldsymbol{\mu}\|^2_{\mathbb{S}'}}, \\
&=& \|q\|^2_{\xLtwo(\Omega)} + \|\boldsymbol{\mu}\|^2_{\mathbb{S}'} \\
&\geq&  \dfrac{1}{\sqrt{C_2}}\, \|\boldsymbol{v}\|_{\mathbf{W}_N}\, \left( \|q\|^2_{\xLtwo(\Omega)} + \|\boldsymbol{\mu}\|^2_{\mathbb{S}'} \right)^{1/2}.
\end{eqnarray*}
The well-posedness of the formulation~\eqref{eq:ddfvmac1}--\eqref{eq:ddfvmac3} follows from the Babu\v{s}ka--Brezzi theorem.
\end{proof}

\medbreak

In order to interpret the decomposed formulation~\eqref{eq:ddfvmac1}--\eqref{eq:ddfvmac3}, we shall need the following lemma.
\begin{lmm}
\label{lambda}
A continuous anti-linear functional $L_W$ on $\mathbf{W}_N$ vanishes on $\mathbf{X}_N(\underline{\boldsymbol{K}};\Omega)$ if, and only if, it is of the form~\eqref{elmu} for some $\boldsymbol{\mu} \in (\mathbb{S}_{\Sigma}^W)' $. More specifically, there exists a unique pair $(\mu_n , \boldsymbol{\mu}_{\top}) \in \xHn{{1/2}}(\Sigma) \times \mathbf{TC}(\Sigma)$ such that:
\begin{eqnarray}
\label{LX}
L_W(\boldsymbol{F}) = \sum_{i \bigtriangleup j} \int_{\Sigma_{i,j}} \left\{ \mu_n \, [\overline{\underline{\boldsymbol{K}} \boldsymbol{F} \cdot  \boldsymbol{n}}]_{\Sigma_{ij}} + \boldsymbol{\mu}_{\top} \cdot [\overline{\boldsymbol{F} \times \boldsymbol{n}}]_{\Sigma_{ij}} \right\}\, \xdif\sigma .
\end{eqnarray}
\end{lmm}
\begin{proof}
Let $L^i_W \in (\mathbf{W}_N^i)'$, then Hahn--Banach and Riesz's theorems (\cf~\cite{Brez83}, Thm~VIII.13) show that there exist $\boldsymbol{g}^i_0 \in \mathbf{L}^2(\Omega_i)$, $g^i_1 \in \xLtwo(\Omega_i)$ and $\boldsymbol{g}^i_2 \in \mathbf{L}^2(\Omega_i)$ such that:
\begin{equation*}
\forall \boldsymbol{F}_i \in \mathbf{W}_N^i,\quad L_W^i(\boldsymbol{F}) = \int_{\Omega_i} \left( \boldsymbol{g}^i_0 \cdot \overline{\boldsymbol{F}} + g^i_1\, \overline{\dive \underline{\boldsymbol{K}} \boldsymbol{F}_i} + \boldsymbol{g}^i_{2} \cdot \overline{\curl \boldsymbol{F}_i} \right)\, \xdif\Omega . 
\end{equation*}
Since $\mathbf{W}_N = \bigoplus_{i} \mathbf{W}_N^i$, we have $(\mathbf{W}_N)' = \bigoplus_{i} (\mathbf{W}_N^i)'$.
It follows that any anti-linear form on $\mathbf{W}_N$ can be written as:
\begin{equation*}
L_W(\boldsymbol{F}) =\sum_{i=1}^{N_d}  \int_{\Omega_i} \left( \boldsymbol{g}_0 \cdot \overline{\boldsymbol{F}} + g_1\, \overline{\dive \underline{\boldsymbol{K}} \boldsymbol{F}} + \boldsymbol{g}_{2} \cdot \overline{\curl \boldsymbol{F}} \right)\, \xdif\Omega ,
\end{equation*}
with $\boldsymbol{g}_0 \in \mathbf{L}^2(\Omega)$, $g_1 \in \xLtwo(\Omega)$ and $\boldsymbol{g}_2 \in \mathbf{L}^2(\Omega)$. 
We perform a Helmholtz decomposition of $\boldsymbol{g}_0$ in $\underline{\boldsymbol{K}}^*$-gradient and solenoidal parts (Remark~\ref{KHgradsole}): 
\begin{equation*}
\boldsymbol{g}_0 = \underline{\boldsymbol{K}}^* \grad \psi + \boldsymbol{g}_T := \boldsymbol{g}_L + \boldsymbol{g}_T, \quad\text{with}\quad
\psi \in \xHone_0(\Omega) \quad\text{and}\quad \dive \boldsymbol{g}_T = 0.
\end{equation*}

\medbreak

Assume that $L_W$ vanishes on $\mathbf{X}_N(\underline{\boldsymbol{K}};\Omega)$. Let $\boldsymbol{F} \in \boldsymbol{\mathcal{D}}(\Omega)$, and consider its Helmholtz decomposition into gradient and $\underline{\boldsymbol{K}}$-solenoidal parts:
\begin{equation*}
\boldsymbol{F} = \grad \phi + \boldsymbol{F}_T  := \boldsymbol{F}_L  + \boldsymbol{F}_T, \quad\text{with}\quad
\phi \in \xHone_0(\Omega) \quad\text{and}\quad \dive \underline{\boldsymbol{K}} \boldsymbol{F}_T = 0.
\end{equation*}
Using~\eqref{Green10}, one immediately checks that $\left( \boldsymbol{g}_L \mid \boldsymbol{F}_T \right)_{\Omega} = 0$ and  $\left( \boldsymbol{g}_T \mid \boldsymbol{F}_L \right)_{\Omega} = 0$.
Furthermore, $\boldsymbol{F}_L$ and $\boldsymbol{F}_T$ belong to $\mathbf{X}_N(\underline{\boldsymbol{K}};\Omega)$, by Remark \ref{contXhelm}. Since $L_W$ vanishes on $\mathbf{X}_N(\underline{\boldsymbol{K}};\Omega)$, we deduce:
\begin{eqnarray}\nonumber
0=L_W(\boldsymbol{F}_L)=\int_{\Omega} \biggl( \boldsymbol{g}_L \cdot \overline{\boldsymbol{F}_L} + \underbrace{\boldsymbol{g}_T \cdot \overline{\boldsymbol{F}_L} }_{=0} \mbox{} + g_1\, \overline{\dive \underline{\boldsymbol{K}} \boldsymbol{F}_L} + \underbrace{\boldsymbol{g}_2 \cdot \overline{\curl \boldsymbol{F}_L} }_{=0} \biggr)\, \xdif\Omega.
\end{eqnarray}
By adding $\int_{\Omega} \left( \boldsymbol{g}_L \cdot \overline{\boldsymbol{F}_T} + g_1\, \overline{\dive \underline{\boldsymbol{K}} \boldsymbol{F}_T} \right)\, \xdif\Omega=0$, we have
\begin{eqnarray}\label{LXFL}
\int_{\Omega} \left( \boldsymbol{g}_L \cdot \overline{\boldsymbol{F}} + g_1\, \overline{\dive \underline{\boldsymbol{K}} \boldsymbol{F}} \right)\, \xdif\Omega = 0, \quad \forall \boldsymbol{F} \in  \boldsymbol{\mathcal{D}}(\Omega) .
\end{eqnarray}
This yields $\underline{\boldsymbol{K}}^* \grad g_1 = \boldsymbol{g}_L$ in $\boldsymbol{\mathcal{D}}'(\Omega)$. As $\boldsymbol{g}_L \in \mathbf{L}^2(\Omega)$ and $\left(\underline{\boldsymbol{K}}^*\right)^{-1} \in \xLinfty(\Omega; \mathcal{M}_3(\xC))$ by Lemma~\Rref{zetaeta}, we infer $g_1 \in \xHone(\Omega)$. Furthermore, $\grad g_1 = \grad \psi$ in~$\Omega$; as $\Omega$ is connected, this gives $g_1 = \psi + C_1$, for some constant~$C_1$. In particular, $g_1=C_1$ on the boundary~$\Gamma$.

\medbreak

Let $w \in \xHone_0(\Omega)$ such that $\Delta_{\underline{\boldsymbol{K}}}w \in \mathbf{L}^2(\Omega)$; then $\grad w \in \mathbf{X}_N(\underline{\boldsymbol{K}};\Omega)$, and we have
\begin{eqnarray}\nonumber
0=L_W(\grad w) &=& \int_{\Omega} \left( \boldsymbol{g}_L \cdot \overline{\grad w} + g_1\, \overline{\Delta_{\underline{\boldsymbol{K}}} w} \right)\, \xdif\Omega \\
&\stackrel{\eqref{Green1}}{=}& \int_{\Omega} \underbrace{(\boldsymbol{g}_L \cdot \overline{\grad w} - \underline{\boldsymbol{K}}^* \grad g_1 \cdot \overline{\grad w} )}_{=0} \xdif\Omega \nonumber \\
& & + \langle g_1 , \underline{\boldsymbol{K}} \grad w \cdot \boldsymbol{n} \rangle_{\xHn{{-1/2}}(\Gamma)} \nonumber
\end{eqnarray}
Taking $w$ such that $\langle 1 , \underline{\boldsymbol{K}} \grad w \cdot \boldsymbol{n} \rangle_{\xHn{{-1/2}}(\Gamma)} \ne 0$, one deduces $C_1 = \left. g_1 \right|_{\Gamma}=0$, \ie, $g_1 \in \xHone_0(\Omega)$. On the other hand, we have
\begin{eqnarray}\nonumber
0=L_W(\boldsymbol{F}_T)=\int_{\Omega} \biggl( \underbrace{\boldsymbol{g}_L \cdot \overline{\boldsymbol{F}_T}}_{=0} \mbox{} + \boldsymbol{g}_T \cdot \overline{\boldsymbol{F}_T} + \underbrace{g_1\, \overline{\dive \underline{\boldsymbol{K}} \boldsymbol{F}_T} }_{=0} \mbox{} + \boldsymbol{g}_2 \cdot \overline{\curl \boldsymbol{F}_T} \biggr)\, \xdif\Omega.
\end{eqnarray}
By adding $\int_{\Omega} \left(\boldsymbol{g}_T \cdot \overline{\boldsymbol{F}_L} + \boldsymbol{g}_2 \cdot \overline{\curl \boldsymbol{F}_L} \right)\, \xdif\Omega = 0$, we have
\begin{eqnarray}\label{LXFT}
\int_{\Omega} \left( \boldsymbol{g}_T \cdot \overline{\boldsymbol{F}} + \boldsymbol{g}_2 \cdot \overline{\curl \boldsymbol{F}} \right)\, \xdif\Omega=0, \quad \forall \boldsymbol{F} \in \boldsymbol{\mathcal{D}}(\Omega).
\end{eqnarray}
Thus, we have $\curl \boldsymbol{g}_2 = - \boldsymbol{g}_T$ in $\boldsymbol{\mathcal{D}}'(\Omega)$ and therefore then $\boldsymbol{g}_2 \in \mathbf{H}(\curl;\Omega)$.
Adding \eqref{LXFL} and \eqref{LXFT}, the form $L_W$ is equal to:
\begin{equation*}
L_W(\boldsymbol{F})= \sum_{i=1}^{N_d} \int_{\Omega_i} \left( g_1\, \overline{\dive \underline{\boldsymbol{K}} \boldsymbol{F}} + \underline{\boldsymbol{K}}^* \grad g_1 \cdot \overline{\boldsymbol{F}} + \boldsymbol{g}_2 \cdot \overline{\curl \boldsymbol{F}} - \curl \boldsymbol{g}_2 \cdot \overline{\boldsymbol{F}} \right)\, \xdif\Omega  .
\end{equation*}
with $g_1 \in \xHone_0(\Omega)$ and $\boldsymbol{g}_2 \in \mathbf{H}(\curl;\Omega)$.
Using the Green formulas \eqref{Green1} and~\eqref{Green2} in each~$\Omega_i$, we deduce
\begin{eqnarray}\nonumber
L_W(\boldsymbol{F})= \sum_{i=1}^{N_d} \int_{\partial \Omega_i} \left( g_1 \left(\overline{\underline{\boldsymbol{K}} \boldsymbol{F} \cdot \boldsymbol{n}} \right)+ \boldsymbol{g}_{2\top} \cdot \left(\overline{\boldsymbol{F} \times \boldsymbol{n}} \right)\right)\, \xdif\sigma.
\end{eqnarray}
Each integral is understood as a sum of two duality products: the first between $\xHn{{1/2}}(\partial\Omega_i)$ and~$\xHn{{-1/2}}(\partial\Omega_i)$; the second between $\mathbf{TC}(\partial\Omega_i)$ and~$\mathbf{TT}(\partial\Omega_i)$.  
On an exterior boundary $\Gamma_i$, one has $g_1=0$ and $\boldsymbol{F} \times \boldsymbol{n} = 0$; on an interface $\Sigma_{i,j}$, the sum of the contributions of $\Omega_i$ and~$\Omega_j$ amounts to a jump, as noted in~\eqref{dual-Sigma}. Finally, we arrive at~\eqref{LX}, where $\mu_n\in \xHn{{1/2}}(\Sigma)$~is the trace of $g_1$ on~$\Sigma$, and $\boldsymbol{\mu}_{\top} \in \mathbf{TC}(\Sigma)$ is the tangential component of $\boldsymbol{g}_2$ on~$\Sigma$. 
These characterisations allow one to consider their restriction to each interface: $\mu_n^{i,j}\in \xHn{{1/2}}(\Sigma_{i,j})$ and $\boldsymbol{\mu}_{\top}^{i,j} \in \mathbf{TC}(\Sigma_{i,j})$ on each~$\Sigma_{i,j}$. Of course, restrictions to neighbouring interfaces satisfy suitable compatibility conditions.

\medbreak

To prove uniqueness, it is enough to show that $L_W(\boldsymbol{F})=0,\ \forall \boldsymbol{F} \in \mathbf{W}_N$, implies $\mu_n = 0$ and $\boldsymbol{\mu}_{\top} = \boldsymbol{0}$. 
First, take $g \in \xHn{{-1/2}}(\Sigma)$, and introduce $\phi \in \xHone_0(\Omega)$ solution to the following variational formulation, with the form~$\mathfrak{a}$ from~\eqref{fv-divKgrad}:
\begin{equation*}
\mathfrak{a}(\phi,\psi) = \langle g , \psi_{|\Sigma} \rangle_{\xHn{{1/2}}(\Sigma)} \, \quad \forall \psi \in \xHone_0(\Omega),
\end{equation*}
which is well-posed as in Lemma~\Rref{divKgrad}.
Performing an integration by parts in each~$\Omega_i$ and adding as before, we see that $\phi$ satisfies:
\begin{equation*} 
-\Delta_{\underline{\boldsymbol{K}}}\phi =  0   \text{ in each } \Omega_i \,, \quad 
[ \phi ]_{\Sigma} = 0 \text{ and } \left[ \underline{\boldsymbol{K}} \grad \phi \cdot \boldsymbol{n} \right]_{\Sigma} = g \text{ on }  \Sigma .
\end{equation*} 
Setting $\boldsymbol{F}=\grad \phi$, we have $\boldsymbol{F} \in \mathbf{W}_N$, $[\underline{\boldsymbol{K}} \boldsymbol{F} \cdot \boldsymbol{n}]_{\Sigma}=g$ and $[ \boldsymbol{F} \times \boldsymbol{n}]_{\Sigma}=0$. So:
\begin{equation*}
0 = L_W (\boldsymbol{F}) = \int_{\Sigma} \mu_{n}\, \left[ \overline{\underline{\boldsymbol{K}} \boldsymbol{F} \cdot \boldsymbol{n}} \right]_{\Sigma}\, \xdif\sigma =  \left\langle \mu_{n}, g \right\rangle_{\xHn{{-1/2}}(\Sigma)}.
\end{equation*}
As $g$ is arbitrary, one deduces $\mu_{n}=0$ in~$\xHn{{1/2}}(\Sigma)$. In particular, taking $g$ supported on one interface~$\Sigma_{i,j}$, one finds $\mu_{n}^{i,j} =0$ in~$\xHn{{1/2}}(\Sigma_{i,j})$.

\medbreak

For the tangential part, take $\boldsymbol{\varphi} \in \mathbf{S}_{\Sigma}^{V}$. By definition, there exists $\boldsymbol{v} \in \mathbf{V}_0$ such that $[\boldsymbol{v} \times \boldsymbol{n}]_\Sigma = \boldsymbol{\varphi}$. In each subdomain~$\Omega_i$, introduce 
$$\phi_i \in \xHone_0(\Omega_i) \text{ solution to: } -\Delta_{\underline{\boldsymbol{K}}}\phi_i = \dive (\underline{\boldsymbol{K}} \boldsymbol{v}_i) \in \xHn{{-1}}(\Omega_i), \quad\text{and}\quad \boldsymbol{F}_i = \boldsymbol{v}_i + \grad \phi_i.$$ 
There holds $\boldsymbol{F}_i \in \mathbf{H}(\curl;\Omega_i) \cap \mathbf{H}(\dive{\underline{\boldsymbol{K}}},\Omega_i)$ and $\boldsymbol{F}_i \times \boldsymbol{n}_i = \boldsymbol{v}_i \times \boldsymbol{n}_i$ on $\partial\Omega_i$. Therefore, the global field $\boldsymbol{F} = \left\{ \boldsymbol{F}_i \right\}_{i=1,\ldots,N_d}$ satisfies $\boldsymbol{F} \in \mathbf{W}_N$ and $[\boldsymbol{F} \times \boldsymbol{n}]_\Sigma = [\boldsymbol{v} \times \boldsymbol{n}]_\Sigma = \boldsymbol{\varphi}$.
As $L_W(\boldsymbol{F})=0$, this implies
$$\langle \boldsymbol{\mu}_{\top} , \boldsymbol{\varphi} \rangle_{\mathbf{S}_{\Sigma}^{V}}=0, \quad \forall \boldsymbol{\varphi} \in \mathbf{S}^{V}_{\Sigma}, \quad\text{\ie,}\quad \boldsymbol{\mu}_{\top} =  \boldsymbol{0} \quad \textrm{in} \quad (\mathbf{S}^{V}_{\Sigma})'.$$
In particular, taking $\boldsymbol{\varphi} \in \widetilde{\mathbf{TT}}(\Sigma_{i,j})$, its extension by~$0$ to~$\Sigma$ belongs to~$\mathbf{S}_{\Sigma}^{V}$ by Lemma~\Rref{TNtilde-SSigmaV}, and we infer $\boldsymbol{\mu}_{\top}^{i,j} =  \boldsymbol{0}$ in $\widetilde{\mathbf{TT}}(\Sigma_{i,j})' = \mathbf{TC}(\Sigma_{i,j})$.
\end{proof}

\medbreak

\begin{thrm}
The decomposed formulation~\eqref{eq:ddfvmac1}--\eqref{eq:ddfvmac3} and the original mixed augmented formulation~\eqref{eq:fvma1}--\eqref{eq:fvma2} are equivalent: $(\boldsymbol{E},p,\boldsymbol{\lambda}) $ is solution to~\eqref{eq:ddfvmac1}--\eqref{eq:ddfvmac3} iff $(\boldsymbol{E},p) $ is solution to~\eqref{eq:fvma1}--\eqref{eq:fvma2}, and
\begin{eqnarray}
\lambda_n = 0,\qquad  \boldsymbol{\lambda}_{\top} = (\curl \boldsymbol{E})_{\top|\Sigma} \,.
\end{eqnarray}
\end{thrm}
\begin{proof}
Let $(\boldsymbol{E},p,\boldsymbol{\lambda}) $ be the solution to~\eqref{eq:ddfvmac1}--\eqref{eq:ddfvmac3}. {F}rom~\eqref{eq:ddfvmac3}, we have the jump conditions~\eqref{saut01}, and $\boldsymbol{E} \in \mathbf{X}_N(\underline{\boldsymbol{K}};\Omega)$. Taking a test function $\boldsymbol{F} \in \mathbf{X}_N(\underline{\boldsymbol{K}};\Omega)$, the term $\langle \boldsymbol{\lambda} , [[\boldsymbol{F} ]]_\Sigma \rangle_{\mathbb{S}_{\Sigma}^W}$ vanishes in~\eqref{eq:ddfvmac1}, which gives~\eqref{eq:fvma1}. Then~\eqref{eq:ddfvmac2} is identical to~\eqref{eq:fvma2}. This means that $(\boldsymbol{E},p) \in \mathbf{X}_N(\underline{\boldsymbol{K}};\Omega) \times \xLtwo(\Omega)$ coincides with the unique solution to~\eqref{eq:fvma1}--\eqref{eq:fvma2}.

\medbreak

Conversely, let $(\boldsymbol{E},p)$ be the solution to~\eqref{eq:fvma1}--\eqref{eq:fvma2}. As $\boldsymbol{E} \in \mathbf{X}_N(\underline{\boldsymbol{K}};\Omega)$, we have automatically $ [\boldsymbol{E} \times \boldsymbol{n}]_{\Sigma} = 0$ and $[\underline{\boldsymbol{K}} \boldsymbol{E} \cdot \boldsymbol{n}]_{\Sigma} = 0$, which implies~\eqref{eq:ddfvmac3}. As for~\eqref{eq:ddfvmac2}, we have
\begin{eqnarray*}
\sum_{i} (\dive \boldsymbol{E}_i \mid q_i) = (\dive \boldsymbol{E} \mid q) = (g \mid q) = \sum_{i}(g_i \mid q_i).
\end{eqnarray*}
Define the continuous anti-linear form $L_W$ on $\mathbf{W}_N$ :
\begin{eqnarray*}
L_W :\quad \boldsymbol{F}  \longmapsto  \sum_i \left( -a_{i,s} ( \boldsymbol{E}_i , \boldsymbol{F}_i)-b_i( \boldsymbol{F}_i , p_i) + L_i( \boldsymbol{F}_i) \right),
\end{eqnarray*}
which vanishes on $\mathbf{X}_N(\underline{\boldsymbol{K}};\Omega)$. By Lemma~\Rref{lambda}, there exists a unique $\boldsymbol{\lambda} \in (\mathbb{S}_{\Sigma}^W)' $ such that
\begin{eqnarray}\nonumber
L_W(\boldsymbol{F}) =  \int_{\Sigma} \left\{ \lambda_n\, [\overline{\underline{\boldsymbol{K}} \boldsymbol{F} \cdot  \boldsymbol{n}}]_{\Sigma} + \boldsymbol{\lambda}_{\top} \cdot [\overline{\boldsymbol{F} \times \boldsymbol{n}}]_{\Sigma} \right\}\, \xdif\sigma .
\end{eqnarray}
So, Eq.~\eqref{eq:ddfvmac1} is verified. On the other hand, we have remarked that the solution to~\eqref{eq:fvma1}--\eqref{eq:fvma2} satisfies $\dive (\underline{\boldsymbol{K}} \boldsymbol{E}) = g$ and $p=0$; thus the strong form of~\eqref{eq:fvma1} becomes:
\begin{equation*}
\curl\curl \boldsymbol{E} - \tfrac{\omega^2}{c^2} \underline{\boldsymbol{K}} \boldsymbol{E} = \boldsymbol{f} \quad \text{in } \boldsymbol{\mathcal{D}}'(\Omega).
\end{equation*}
As a consequence, $\curl\boldsymbol{E} \in \mathbf{H}(\curl;\Omega)$. 
Starting again from~\eqref{eq:ddfvmac1}, using the Green formulas \eqref{Green2},~\eqref{Green1} in each~$\Omega_i$, and taking the above equalities into account, one obtains:
\begin{equation*}
\langle -\curl \boldsymbol{E} , [\boldsymbol{F} \times \boldsymbol{n}]_{\Sigma} \rangle_{\mathbf{S}_{\Sigma}^{V}} 
+ \langle \lambda_n , [ \underline{\boldsymbol{K}} \boldsymbol{F} \cdot \boldsymbol{n} ]_{\Sigma} \rangle_{\xHn{{-1/2}}(\Sigma)} 
+ \langle \boldsymbol{\lambda}_{\top} , [\boldsymbol{F} \times \boldsymbol{n}]_{\Sigma} \rangle_{\mathbf{S}_{\Sigma}^{V}}  
= 0.
\end{equation*}
Thus we get the expressions of $\lambda_n$ and $\boldsymbol{\lambda}_{\top}$.
\end{proof}

\begin{acknowledgement} 
Acknowledgement: The authors wish to thank the anonymous referees for their useful remarks and suggestions.
\end{acknowledgement}


\end{document}